\documentclass[a4paper,10pt]{article}
\usepackage[a4paper, total={7in, 10in}]{geometry}
\usepackage[utf8x]{inputenc}
\usepackage{comment}
\usepackage{amssymb}
\usepackage{amsthm}
\usepackage{hyperref}
\theoremstyle{plain}
\newtheorem{thm}{Theorem}[section]
\newtheorem{lem}[thm]{Lemma}

\newtheorem{cor}{Corollary}

\theoremstyle{defn}
\newtheorem{defn}{Definition}[section]

\theoremstyle{rmk}
\newtheorem{rmk}{Remark}

\newtheorem{pf}{Proof}
\usepackage{amsmath}
\usepackage{enumitem}
\title{Two Timescale Stochastic Approximation with Controlled Markov noise and Off-policy Temporal Difference Learning}

\author{Prasenjit Karmakar and Shalabh Bhatnagar}

\begin{document}

\maketitle

\begin{abstract}
We present for the first time an asymptotic convergence analysis of two time-scale stochastic approximation driven by
`controlled' Markov noise. In particular, both the faster and slower recursions have non-additive controlled Markov noise components in 
addition to martingale difference noise. We analyze the asymptotic behavior of our framework
by relating it
to limiting differential inclusions in both time-scales that are defined in terms 
of the ergodic occupation measures associated with the controlled 
Markov processes. 
Finally, we present a solution to the off-policy convergence problem for temporal difference 
learning with linear function approximation, using 
our results. 
\end{abstract}



\maketitle

%

\section{Introduction}

Stochastic approximation algorithms are sequential non-parametric methods for finding a zero or minimum of a function in the situation where
only the noisy observations of the function values are available. Two time-scale stochastic approximation algorithms represent
one of the most general subclasses of stochastic approximation methods. These algorithms consist of two coupled recursions which 
are updated with different (one is considerably smaller than the other) step sizes which in turn facilitate convergence for such algorithms.

Two time-scale stochastic approximation algorithms \cite{borkartt} have successfully been applied to several complex problems arising in the areas of
reinforcement learning, signal processing and admission control in communication networks. There are many reinforcement learning applications
(precisely those where parameterization of value function is implemented) where non-additive Markov noise is present in one or both iterates thus 
requiring the current two time-scale framework to be extended to
include Markov noise (for example, in \cite[p.~5]{sutton2} it is mentioned that in order to generalize the analysis 
to Markov noise, the theory of two time-scale stochastic approximation needs to include the latter).

Here we present a more general framework of two time-scale stochastic approximation
with ``controlled'' Markov noise, i.e., the noise is 
not simply Markov; rather it is driven by the 
iterates and an additional control process as well. 
We analyze the asymptotic behaviour of our framework by relating it to limiting differential inclusions in both timescales that 
are defined in terms of the ergodic occupation measures associated with the controlled Markov processes. Next, using these 
results 
for the special
case of our framework where the random processes are irreducible Markov chains, 
we present a solution to the off-policy convergence problem
for temporal difference learning with linear function approximation. While the 
off-policy convergence problem for reinforcement learning (RL) with linear 
function approximation has
been one of the most interesting problems, there are very few solutions available in the
current literature. One such work \cite{yu} shows 
the convergence of the least squares temporal difference learning algorithm with 
eligibility traces (LSTD($\lambda$)) as well as the TD($\lambda$) algorithm. 
While the LSTD methods are not feasible when 
the dimension of the feature vector is large, off-policy TD($\lambda$) is shown to 
converge only when the eligibility function $\lambda \in [0,1]$ is very close to 1. Another recent work \cite{yu_new}
proves weak convergence of several emphatic temporal difference learning algorithms which is also designed to 
solve the off-policy convergence problem.
In \cite{sutton1,sutton,maeith} the 
gradient temporal difference learning (GTD) algorithms were proposed to solve this problem. However, the authors 
make the assumption that the data
is available in the ``off-policy'' setting (i.e. the off-policy issue is incorporated into the data rather than in the algorithm) 
whereas, in reality, one has only the ``on-policy'' Markov trajectory 
corresponding to a given behaviour policy and we are interested in designing an online learning algorithm. We use one 
of the algorithms from \cite{maeith} called TDC with ``importance-weighting'' which takes the ``on-policy'' data as 
input and show its convergence using the results we develop. 
Our convergence analysis can also be extended for the same algorithm with 
eligibility traces for a sufficiently large range of values of $\lambda$. 
Our results can be used to provide a convergence analysis for reinforcement learning 
algorithms such as those in \cite{mannor} for which convergence proofs have not been provided. 
 
To the best of our knowledge there are related works  such as \cite{Tadic, konda, konda_actor, tadic_new} 
where two time-scale stochastic approximation 
algorithms with algorithm iterate dependent non-additive Markov noise is analyzed. In all of them the 
Markov noise in the recursion is handled using the classic Poisson equation based approach of \cite{benveniste, metivier} and 
applied to the asymptotic analysis of many algorithms used in machine learning, system identification, signal 
processing, image analysis and automatic control. However, we show that our method 
also works if there is another additional control process as well and if the
underlying Markov process has non-unique stationary distributions. Further, the mentioned application 
does not 
require strong assumption such as aperiodicity for the underlying Markov chain which 
is a sufficient condition if we use Poisson equation based approach \cite{adam, Tadic}.  Additionally, our assumptions 
are quite different from the assumptions made in the mentioned literature and we give a detailed
 comparison in Section \ref{defn}.

The organization of the paper is as follows: Section \ref{sec_def} formally defines the problem and 
provides background and assumptions. Section \ref{mres} shows the main
results. Section \ref{relax} discusses how one of our assumptions of Section \ref{sec_def} can be relaxed. 
Section \ref{app} presents an application of our results to the off-policy convergence problem for temporal difference 
learning with linear function approximation. 
Finally, we conclude by providing
some future research directions. 
\section{Background, Problem Definition, and Assumptions}
\label{sec_def}

In the following we describe the preliminaries and notation used in our proofs. 
Most of the definitions and notation are from \cite{benaim,borkar1,Aubin}.
\subsection{\textit{Definition and Notation}}

Let $F$ denote a set-valued function mapping each point $\theta \in \mathbb{R}^m$ to a set $F(\theta) \subset \mathbb{R}^m$. $F$ is called
a \textit{Marchaud map} if the following hold: 
\begin{enumerate}[label=(\roman*)]
\item  $F$ is \textit{upper-semicontinuous} in the sense that if $\theta_n \to \theta$ and $w_n \to w$ with $w_n \in F(\theta_n)$ for all $n\geq 1$, then 
$w \in F(\theta)$. In order words, the graph of $F$ defined as $\{(\theta,w):  w \in F(\theta)\}$ is closed. 
\item $F(\theta)$ is a non-empty compact convex subset of $\mathbb{R}^m$ for all $\theta \in \mathbb{R}^m$.
\item $\exists c >0$ such that for all $\theta \in \mathbb{R}^m$,
\begin{equation}
\sup_{z\in F(\theta)} \|z\| \leq c(1+\|\theta\|),\nonumber 
\end{equation}
where $\|.\|$ denotes any norm on $\mathbb{R}^m$. 
\end{enumerate}


\textit{A solution for the differential inclusion (d.i.)} 
\begin{equation}
\label{diffin}
\dot{\theta}(t) \in F(\theta(t)) 
\end{equation}
with initial point $ \theta_0 \in \mathbb{R}^m$ is an absolutely continuous (on compacts) 
mapping  $\theta :  \mathbb{R} \to \mathbb{R}^m$ such that
$\theta(0) =\theta_0$ and 
\begin{equation}
\dot{\theta}(t) \in F(\theta(t))\nonumber 
\end{equation}
for almost every $t \in \mathbb{R}$. 
If $F$ is a Marchaud map, it is well-known that (\ref{diffin}) has solutions 
(possibly non-unique) through every initial point. The differential inclusion  (\ref{diffin}) 
induces a \textit{set-valued dynamical system} $\{\Phi_t\}_{t\in \mathbb{R}}$ defined by
\begin{equation}
\Phi_t(\theta_0) = \{\theta(t) : \theta(\cdot) \mbox{ is a solution to  
(\ref{diffin}) with $\theta(0) =\theta_0$}\}.\nonumber 
\end{equation}  \indent
Consider the autonomous ordinary differential equation (o.d.e.)
\begin{equation}
\label{ode1}
\dot{\theta}(t)=h(\theta(t)), 
\end{equation}
where $h$ is Lipschitz continuous. One can write (\ref{ode1}) in the 
format of (\ref{diffin}) by taking $F(\theta)=\{h(\theta)\}$. 
It is well-known that (\ref{ode1}) is well-posed, i.e., it has a \textit{unique solution} for every initial point.
Hence the set-valued dynamical system induced by the o.d.e. or \textit{flow} is  $\{\Phi_t\}_{t\in \mathbb{R}}$ with 
\begin{equation}
\Phi_t(\theta_0) = \{\theta(t)\},\nonumber 
\end{equation}
where $\theta(\cdot)$ is the solution to  
(\ref{ode1}) with $\theta(0) =\theta_0$.
It is also well-known that $\Phi_t(.)$ is a \textit{continuous function} for all $t \in \mathbb{R}$.
\\ \indent     
A set $A \subset \mathbb{R}^m$ is said to be \textit{invariant} (for $F$) if for all 
$\theta_0\in A$ there exists a solution $\theta(\cdot)$ of (\ref{diffin})
with $\theta(0) = \theta_0$ such that $\theta(\mathbb{R}) \subset A$. 
\\ \indent
Given a set $A \subset \mathbb{R}^m$ and $\theta'',w''\in A$, we write $\theta''\hookrightarrow_A w''$ 
if for every $\epsilon > 0$ and $T>0$ $\exists n \in \mathbb{N}$, solutions $\theta_1(\cdot), \dots, \theta_n(\cdot)$ 
to  (\ref{diffin}) and real numbers $t_1, t_2, \dots, t_n$ greater than 
$T$ such that
\begin{enumerate}[label=(\roman*)]
 \item  $\theta_i(s) \in A$ for all $0 \leq s \leq t_i$ and for all $i=1, \dots, n,$
 \item  $\|\theta_i(t_i) - \theta_{i+1}(0)\| \leq \epsilon$ for all $i=1, \dots, n-1,$
 \item  $\|\theta_1(0) - \theta''\| \leq \epsilon$ and $\|\theta_n(t_n) - w''\| \leq \epsilon.$
\end{enumerate}
The sequence $(\theta_1(\cdot), \dots, \theta_n(\cdot))$ is called an $(\epsilon, T)$ 
chain (in $A$ from $\theta''$ to $w''$) for $F$. A set $A \subset \mathbb{R}^m$
is said to be \textit{internally chain transitive}, provided that $A$ is compact and $\theta'' 
\hookrightarrow_A w''$ for all $\theta'',w''\in A$. It can be 
proved that in the above case, $A$ is an invariant set. 
\\ \indent
A compact invariant set $A$ is called an \textit{attractor} for $\Phi$, provided that there is a neighbourhood $U$ of $A$
(i.e., for the induced topology) with the property that 
$d(\Phi_t(\theta''), A) \to 0$ as $t\to \infty$ \textit{uniformly} in $\theta'' \in U$.
Here $d(X, Y) = \sup_{\theta'' \in X}\inf_{w'' \in Y}\|\theta''-w''\|$ for $X,Y \subset \mathbb{R}^m.$
Such a $U$ is called a \textit{fundamental neighbourhood} of the attractor $A$.
\textit{An attractor of a well-posed o.d.e.} is an attractor for the set-valued dynamical system induced by the o.d.e.    
\\ \indent
The set 
\begin{equation}
\omega_\Phi(\theta'') = \bigcap_{t \geq 0} \overline{\Phi_{[t, \infty)}(\theta'')} \nonumber 
\end{equation}
is called the \textit{$\omega$-limit} set of a point $\theta'' \in \mathbb{R}^m$. 
If $A$ is a set, then 
\begin{equation}
B(A) = \{\theta'' \in \mathbb{R}^m :  \omega_\Phi(\theta'') \subset A\} \nonumber 
\end{equation}
denotes its \textit{basin of attraction}. A \textit{global attractor} for $\Phi$ is an attractor $A$
whose basin of attraction consists of all $\mathbb{R}^m$. 
Then the following lemma will be useful for our proofs, see \cite{benaim} for a proof.
\begin{lem}
\label{ga}
Suppose $\Phi$ has a global attractor $A$. Then every internally chain transitive set lies in $A$.
\end{lem}

We also require another result which will be useful to apply our results to the RL application we mention. Before stating it 
we recall some definitions from Appendix 11.2.3 of \cite{borkar1}:
\\ \indent
A point $\theta^*\in \mathbb{R}^m$ is called \textit{Lyapunov stable} for the o.d.e (\ref{ode1}) if for all $\epsilon >0$, 
there exists a $\delta >0$ such that every trajectory of (\ref{ode1}) initiated in the 
$\delta$-neighbourhood of $\theta^*$ remains in its $\epsilon$-neighbourhood. 
$\theta^*$ is called \textit{globally asymptotically stable} if $\theta^*$ is Lyapunov stable and 
\textit{all} trajectories of the o.d.e. converge to it. 
\begin{lem}
\label{ga2}
Consider the autonomous o.d.e. $\dot{\theta}(t)=h(\theta(t))$ where $h$ is Lipschitz continuous. Let $\theta^*$ be 
globally asymptotically stable. Then $\theta^*$ is the global 
attractor of the o.d.e. 
\end{lem}
\begin{pf}
We refer the readers to Lemma~1 of \cite[Chapter 3]{borkar1} for a proof. 
\end{pf}

We end this subsection with a notation which will be used frequently in the convergence statements in the following 
sections.
\begin{defn}
For function $\theta(.)$ defined on $[0, \infty)$, the notation ``$\theta(t) \to A$ as $t \to \infty$'' 
means that $\cap_{t \geq 0} \overline{\{\theta(s):s \geq t\}} \subset A$.
Similar definition applies for a sequence $\{\theta_n\}$.    
\end{defn}

\subsection{\textit{Problem Definition}}
\label{defn}

Our goal is to perform an asymptotic analysis of the following coupled recursions:  
\begin{eqnarray}
\theta_{n+1} &= \theta_n + a(n)\left[h(\theta_n, w_n, Z^{(1)}_n) + M^{(1)}_{n+1}\right],\label{eqn1}\\
w_{n+1} &= w_n + b(n)\left[g(\theta_n, w_n, Z^{(2)}_n) + M^{(2)}_{n+1}\right],\label{eqn2}
\end{eqnarray}
where $\theta_n \in \mathbb{R}^d, w_n \in \mathbb{R}^k, n\geq 0$ and $\{Z^{(i)}_n\}, \{M^{(i)}_{n}\}, i=1, 2$ 
are random processes that we describe below. 
\\ \indent 
We make the following assumptions:  
\begin{enumerate}[label=\textbf{(A\arabic*)}]
 \item $\{Z^{(i)}_n\}$ takes values in a compact metric space $S^{(i)}, i=1,2$. Additionally, 
the processes $\{Z^{(i)}_n\}, i = 1, 2$ are controlled 
Markov processes that are controlled by three different control processes: the iterate sequences $\{\theta_m\}, \{w_m\}$ and a 
random process $\{A^{(i)}_n\}$ taking values in a compact metric space $U^{(i)}$
respectively with their individual dynamics specified by
\begin{equation}
P(Z^{(i)}_{n+1} \in B^{(i)} |Z^{(i)}_m, A^{(i)}_m, \theta_m, w_m, m\leq n) = \int_{B^{(i)}} p^{(i)}(dy|Z^{(i)}_n, A^{(i)}_n, \theta_n, w_n), n\geq 0, \nonumber 
\end{equation} 
for $B^{(i)}$ Borel in $S^{(i)}, i = 1, 2,$ respectively.  

\begin{rmk}
In this context one should note that 
\cite{benveniste, metivier} require the Markov process to take values in a normed Polish space. 
\end{rmk}

\begin{rmk}
In \cite{borkar} it is assumed that the state space where the controlled Markov Process takes values is Polish. This space is then compactified
using the fact that a Polish space can be homeomorphically embedded into a dense subset of a compact metric space. The vector 
field $h(.,.) :  \mathbb{R}^d \times S \to \mathbb{R}^d$ is considered bounded when 
the first component lies in a compact set. This would, however, require a continuous 
extension of $h': \Bbb R^d \times \phi(S) \to \Bbb R^d$ defined by $h'(x,s') = h(x,\phi^{-1}(s'))$ 
to $\Bbb R^d \times \overline{\phi(S)}$. 
Here $\phi(\cdot)$ is the homeomorphism defined by 
$\phi(s) = (\rho(s, s_1), \rho(s, s_2), \dots) \in [0,1]^{\infty}$, and 
$\{s_i\}$ and $\rho$ is a countable dense subset and metric of the Polish space
 respectively. A sufficient condition for the above  
is $h'$ to be uniformly continuous \cite[Ex:13, p.~99]{Rudin}.  
However, this is hard to verify.  
This is the main motivation for us to take the range of the Markov process
as compact for our problem. However, there are other reasons for taking compact state space which will be clear 
in the proofs of this section and the next. 
\end{rmk}

\item $h :  \mathbb{R}^{d+k} \times S^{(1)} \to \mathbb{R}^d$ is  
jointly continuous as well as Lipschitz in its first two arguments uniformly w.r.t the third. The latter condition means that
\begin{equation}
\forall z^{(1)} \in S^{(1)}, \|h(\theta, w, z^{(1)}) - h(\theta', w', z^{(1)})\| \leq L^{(1)}(\|\theta-\theta'\| + \|w - w'\|).\nonumber
\end{equation}
Same thing is also true for $g$ where the Lipschitz constant is $L^{(2)}$.
Note that the Lipschitz constant $L^{(i)}$ does not depend on $z^{(i)}$ for $i=1,2$.
\begin{rmk}
We later relax the uniformity of the Lipschitz constant w.r.t the Markov process state space by putting 
suitable moment assumptions on the Markov process.  
\end{rmk}

\item $\{M^{(i)}_n\}, i=1, 2$ are martingale difference sequences
w.r.t increasing $\sigma$-fields
\begin{equation}
\mathcal{F}_n = \sigma(\theta_m, w_m, M^{(i)}_{m}, Z^{(i)}_m, m \leq n, i = 1, 2), n \geq 0,\nonumber 
\end{equation}
satisfying 
\begin{equation}
E[\|M^{(i)}_{n+1}\|^2|\mathcal{F}_n] \leq K(1 + \|\theta_n\|^2 + \|w_n\|^2), i = 1, 2,\nonumber 
\end{equation}
for $n \geq 0$ and a given constant $K>0$.
\item The stepsizes $\{a(n)\}, \{b(n)\}$ are positive scalars satisfying
\begin{equation}
\sum_n a(n) = \sum_n b(n) = \infty, \sum_{n}(a(n)^2 + b(n)^2) < \infty, \frac{a(n)}{b(n)} \to 0.\nonumber 
\end{equation}
Moreover, $a(n), b(n), 
n \geq 0$ are non-increasing. 

Before stating the assumption on the transition kernel $p^{(i)}, i=1, 2$ we need
to define the metric in the space of probability measures $\mathcal{P}(S)$. Here we mention the definitions
and main theorems on the spaces of probability measures that we use in our proofs 
(details can be found in Chapter 2 of \cite{borkar2}). 
We denote the metric by $d$ and is defined as
\begin{equation}
d(\mu, \nu) = \sum_{j} 2^{-j}|\int f_j d\mu - \int f_j d\nu|, \mu, \nu \in \mathcal{P}(S),\nonumber 
\end{equation}
where $\{f_j\}$ are countable dense in the unit ball of $C(S)$.    
Then the following are equivalent: 
\begin{enumerate}[label=(\roman*)]
 \item $d(\mu_n, \mu) \to 0,$
 \item for all bounded  $f$ in $C(S)$, 
\begin{equation}
\int_{S} fd\mu_n \to \int_{S} f d\mu,  
\end{equation} 
\item for all $f$ bounded and uniformly continuous, 
\begin{equation}
\int_{S} fd\mu_n \to \int_{S} fd\mu.\nonumber 
\end{equation} 
\end{enumerate}
Hence we see that $d(\mu_n, \mu) \to 0$ iff 
$\int_{S} f_jd\mu_n \to \int_{S} f_jd\mu$ for all $j$. Any such sequence of functions $\{f_j\}$ is called a convergence
determining class in $\mathcal{P}(S)$. 
Sometimes we also denote $d(\mu_n, \mu) \to 0$ using the notation $\mu_n \Rightarrow \mu$.  
\\ \indent
Also, we recall the characterization of relative compactness in $\mathcal{P}(S)$ that relies on the 
definition of tightness. $\mathcal{A}\subset\mathcal{P}(S)$ is 
a tight set if for any $\epsilon >0$, there exists a compact $K_\epsilon \subset S$ such that $\mu(K_\epsilon) > 1-\epsilon$ for all $\mu \in \mathcal{A}$. 
Clearly, if $S$ is compact then any $\mathcal{A}\subset\mathcal{P}(S)$ is tight. 
By Prohorov's theorem, $\mathcal{A}\subset\mathcal{P}(S)$ is relatively compact if and only if it is tight.
\\ \indent
With the above definitions we assume the following:
\item The map $S^{(i)} \times U^{(i)} \times \mathbb{R}^{d+k} \ni (z^{(i)}, a^{(i)}, \theta, w)  
\to p^{(i)}(dy|z^{(i)}, a^{(i)}, \theta,w) \in \mathcal{P}(S^{(i)})$ is continuous. 
\begin{rmk}
\textbf{(A5)} is much simpler than the assumptions on $n$-step 
transition kernel in \cite[Part II,Chap. 2, Theorem 6]{benveniste}.
\end{rmk}

Additionally, unlike \cite[p~140 line 13]{borkar}, we do not require the extra assumption of the 
continuity in the $\theta$ variable of $p(dy|z,a,\theta)$ to be uniform on compacts w.r.t 
the other variables.

For $\theta_n = \theta, w_n = w$ for all $n$ with a fixed deterministic 
$(\theta, w) \in \mathbb{R}^{d+k}$ and under any stationary randomized control $\pi^{(i)}$, 
it follows from Lemma 2.1 and Lemma 3.1 of \cite{borkar}
that 
the time-homogeneous
Markov processes $Z^{(i)}_n, i=1, 2$ have (possibly non-unique) invariant 
distributions $\Psi^{(i)}_{\theta,w,\pi^{(i)}}, i = 1, 2$. 
. 

Now, it is well-known that the ergodic occupation measure defined as 
\begin{equation}
\Psi^{(i)}_{\theta, w, \pi^{(i)}}(dz, da):= \Psi^{(i)}_{\theta,w,\pi^{(i)}}(dz) \pi^{(i)}(z, da) \in \mathcal{P}(S^{(i)} \times U^{(i)}) \nonumber
\end{equation}satisfies the following: 
\begin{equation}
\label{eqn3} 
\int_{S^{(i)}}f^{(i)}(z) \Psi^{(i)}_{\theta, w, \pi^{(i)}}(dz, U^{(i)}) = \int_{S^{(i)}\times U^{(i)}}\int_{S^{(i)}}f^{(i)}(y)p^{(i)}(dy|z,a, \theta, w)\Psi^{(i)}_{\theta, w, \pi^{(i)}}(dz, da)
\end{equation}
for $f^{(i)}:  S^{(i)} \to \mathcal{R} \in C_b(S^{(i)})$.
\end{enumerate}
We denote by $D^{(i)}(\theta,w), i=1,2$ the set of all such ergodic occupation measures for the prescribed $\theta$ and $w$. In the following we prove
some properties of the map $(\theta,w) \to D^{(i)}(\theta,w)$.

\begin{lem}
\label{lemma1}
For all $(\theta,w)$, $D^{(i)}(\theta,w)$ is convex and compact.
\end{lem}
\begin{pf}
The proof trivially follows from \textbf{(A1)}, \textbf{(A5)} and (\ref{eqn3}).
\end{pf}

\begin{lem}
\label{upsem}
The map $(\theta,w) \to D^{(i)}(\theta,w)$ is upper-semi-continuous.
\end{lem}
\begin{pf}
Let $\theta_n \to \theta, w_n \to w$ and $\Psi^{(i)}_n \Rightarrow \Psi^{(i)} \in \mathcal{P}(S^{(i)} \times U^{(i)})$ such that  $\Psi^{(i)}_n \in D^{(i)}(\theta_n, w_n)$. Let 
$g^{(i)}_n(z,a) = \int_{S^{(i)}}f^{(i)}(y)p^{(i)}(dy|z,a, \theta_n, w_n)$ and  $g^{(i)}(z,a) = \int_{S^{(i)}}f^{(i)}(y)p^{(i)}(dy|z,a, \theta, w).$
From (\ref{eqn3}) we get that
\begin{align} 
\int_{S^{(i)}}f^{(i)}(z) \Psi^{(i)}(dz,U^{(i)}) &= \lim_{n\to \infty}\int_{S^{(i)}}f^{(i)}(z) \Psi^{(i)}_n(dz,U^{(i)})\nonumber\\ 
                                        &= \lim_{n\to \infty}\int_{S^{(i)}\times U^{(i)}}\int_{S^{(i)}}f^{(i)}(y)p^{(i)}(dy|z,a, \theta_n, w_n)\Psi^{(i)}_n(dz,da)\nonumber\\        
                                        &= \lim_{n\to \infty}\int_{S^{(i)}\times U^{(i)}}g^{(i)}_n(z,a)\Psi^{(i)}_n(dz,da).\nonumber
\end{align}
Now, $p^{(i)}(dy|z,a, \theta_n, w_n) \Rightarrow p^{(i)}(dy|z,a, \theta, w)$ 
implies $g^{(i)}_n(\cdot, \cdot) \to  g^{(i)}(\cdot, \cdot)$ pointwise. We prove that the 
convergence is indeed uniform. It is enough to prove that this sequence 
of functions is equicontinuous. Then along with pointwise convergence
it will imply uniform convergence on compacts \cite[p.~168, Ex: 16]{Rudin}. This is also 
a place where \textbf{(A1)} is used.  
\\ \indent
Define $g' : S^{(i)} \times U^{(i)} \times \mathbb{R}^{d+k} \to \mathbb{R}$ by $g'(z',a', \theta',w')  =  \int_{S^{(i)}}f^{(i)}(y)p^{(i)}(dy|z,a', \theta', w')$. 
Then $g'$ is continuous. Let $A= S^{(i)} \times U^{(i)} \times (\{\theta_n\} \cup \theta) \times (\{w_n\} \cup w)$.
So, $A$ is compact and $g'|_{A}$ is uniformly continuous. 
This implies that for all $\epsilon >0$, there exists $\delta >0 $ 
such that if $\rho'(s_1, s_2) < \delta, \mu'(a_1, a_2) < \delta, \|\theta_1-\theta_2\| < \delta,   
\|w_1-w_2\| < \delta,$ then $|g'(s_1, a_1, \theta_1, w_1) - g'(s_2, a_2, \theta_2, w_2)|< \epsilon$  where $s_1, s_2 \in S^{(i)}, 
a_1, a_2 \in U^{(i)},
 \theta_1, \theta_2 \in (\{\theta_n\} \cup \theta), w_1, w_2 \in(\{w_n\} \cup w)$
and $\rho'$ and $\mu'$ denote the metrics in $S^{(i)}$ and $U^{(i)}$ respectively. 
Now use this same $\delta$ for the $\{g^{(i)}_n(\cdot, \cdot)\}$ to get for all $n$ the following for $\rho'(z_1, z_2) < \delta, 
\mu'(a_1, a_2) < \delta$: 
\begin{align}
|g^{(i)}_n(z_1,a_1) - g^{(i)}_n(z_2,a_2)| = |g'(z_1,a_1, \theta_n,w_n) - g'(z_2,a_2, \theta_n,w_n)| < \epsilon.\nonumber 
\end{align}
Hence $\{g^{(i)}_n(\cdot, \cdot)\}$ is equicontinuous. For large $n$, $\sup_{(z,a) \in S^{(i)} \times U^{(i)}}|g^{(i)}_n(z,a) - g^{(i)}(z,a)| < \epsilon/2$ 
because of uniform convergence of $\{g^{(i)}_n(\cdot, \cdot)\}$, hence $\int_{S^{(i)}\times U^{(i)}}|g^{(i)}_n(z,a) - g^{(i)}(z,a)|\Psi^{(i)}_n(dz,da) < \epsilon/2$. 
Now (for $n$ large),
\begin{align}
\label{limit} 
&|\int_{S^{(i)}\times U^{(i)}}g^{(i)}_n(z,a)\Psi^{(i)}_n(dz,da) - \int_{S^{(i)}\times U^{(i)}}g^{(i)}(z,a)\Psi^{(i)}(dz,da)|\nonumber\\ 
&= |\int_{S^{(i)}\times U^{(i)}}[g^{(i)}_n(z,a) - g^{(i)}(z,a)] \Psi^{(i)}_n(dz,da) + \int_{S^{(i)}\times U^{(i)}}g^{(i)}(z,a)\Psi^{(i)}_n(dz,da) - \int_{S^{(i)}\times U^{(i)}}g^{(i)}(z,a)\Psi^{(i)}(dz,da)|\nonumber\\  
&< \epsilon/2 + |\int_{S^{(i)}\times U^{(i)}}g^{(i)}(z,a)\Psi^{(i)}_n(dz,da) - \int_{S^{(i)}\times U^{(i)}}g^{(i)}(z,a)\Psi^{(i)}(dz,da)|\nonumber\\
&< \epsilon.
\end{align}
 The last inequality follows the fact that $\Psi^{(i)}_n \Rightarrow \Psi^{(i)}$. Hence from  (\ref{limit}) we get,
\begin{align} 
\int_{S^{(i)}}f^{(i)}(z) \Psi^{(i)}(dz, U^{(i)}) = \int_{S^{(i)}\times U^{(i)}}\int_{S^{(i)}}f^{(i)}(y)p^{(i)}(dy|z,a, \theta, w)\Psi^{(i)}(dz,da)\nonumber
\end{align} 
proving that the map is upper-semi-continuous.   
\end{pf}

Define $\tilde{g}(\theta, w, \nu) = \int g(\theta,w,z)\nu(dz, U^{(2)})$ for $\nu \in P(S^{(2)}\times U^{(2)})$ 
and $\hat{g}(\theta,w) = \{\tilde{g}(\theta, w, \nu):  \nu \in D^{(2)}(\theta, w)\}.$
\begin{lem}
\label{march}
$\hat{g}(\cdot, \cdot)$ is a Marchaud map. 
\end{lem}
\begin{pf}
\begin{enumerate}[label=(\roman*)]
  \item Convexity and compactness follow trivially from the same for the map $(\theta,w) \to D^{(2)}(\theta, w)$.
  \item \begin{align}
         &\|\tilde{g}(\theta, w, \nu)\|\nonumber\\
         &= \|\int g(\theta,w,z)\nu(dz,U^{(2)})\|\nonumber\\
         &\leq \int\|g(\theta,w,z)\|\nu(dz,U^{(2)})\nonumber\\
         &\leq \int L^{(2)}(\|\theta\|+ \|w\|+\|g(0,0,z)\|)\nu(dz,U^{(2)})\nonumber\\
         &\leq \max (L^{(2)}, L^{(2)}\int \|g(0,0,z)\| \nu(dz,U^{(2)}))(1+\|\theta\|+\|w\|).\nonumber
        \end{align} Clearly, $K(\theta) = \max (L^{(2)}, L^{(2)}\int \|g(0,0,z)\| \nu(dz,U^{(2)})) > 0$. The 
above is true for all $\tilde{g}(\theta,w, \nu) \in \hat{g}(\theta,w), \nu \in D^{(2)}(\theta, w)$.
   \item Let $(\theta_n,w_n) \to (\theta,w), \tilde{g}(\theta_n, w_n, \nu_n) \to m, \nu_n \in D^{(2)}(\theta_n, w_n)$. 
Now, $\{\nu_n\}$ is tight, hence has a convergent
sub-sequence $\{\nu_{n_k}\}$ with $\nu$ being the limit. Then using the arguments similar to the proof of Lemma \ref{upsem} 
one can show that  $m=\tilde{g}(\theta, w, \nu)$
whereas $\nu \in D^{(2)}(\theta, w)$ follows
directly from the upper-semi-continuity of the map $(\theta,w) \to D^{(2)}(\theta,w)$ for all $\theta$. 
\end{enumerate}
\end{pf}
Note that the map $\hat{h}(\cdot,\cdot)$ can be defined similarly and can be shown to be a 
Marchaud map using the exact same technique.

\subsection{\textit{Other assumptions needed for two time-scale convergence analysis}}
\label{assump}

We now list the other assumptions required for two time-scale convergence analysis:  
\begin{enumerate}[label=\textbf{(A\arabic*)}]
\setcounter{enumi}{5}
\item for all $\theta \in \mathbb{R}^d$, the differential inclusion
\begin{equation}
\label{fast}
\dot{w}(t) \in \hat{g}(\theta,w(t)) 
\end{equation}
has a singleton global attractor $\lambda(\theta)$ 
where $\lambda :  \mathbb{R}^d \to \mathbb{R}^k$ is a Lipschitz map with constant $K$.
Additionally,  there exists a continuous function $V: \mathbb{R}^{d+k} \to [0,\infty)$ satisfying the hypothesis of 
Corollary 3.28 of \cite{benaim} with $\Lambda = \{(\theta,\lambda(\theta)):\theta \in \mathbb{R}^d\}$.
This is the most important
assumption as it links the fast and slow iterates.

\item Stability of the iterates: $\sup_n(\|\theta_n\| + \|w_n\|) < \infty$ a.s. 
\end{enumerate}

Let $\bar{\theta}(.), t\geq 0$ be the continuous, piecewise linear 
trajectory defined by $\bar{\theta}(t(n))=\theta_n, n\geq 0$, with 
linear interpolation on each interval $[t(n), t(n+1))$, i.e.,  
\begin{equation}
\bar{\theta}(t) = \theta_n + (\theta_{n+1} - \theta_n)\frac{t-t(n)}{t(n+1)-t(n)}, t \in [t(n), t(n+1)).\nonumber 
\end{equation}
The following theorem is our main result:
\begin{thm}[Slower timescale result]Under assumptions \textbf{(A1)-(A7)}, 
\label{thm}
\begin{equation}
(\theta_n, w_n) \to \cup_{\theta^* \in A_0}(\theta^*, \lambda(\theta^*)) \mbox{a.s. as $n \to \infty$.},\nonumber 
\end{equation}
\end{thm}
where $A_0 = \cap_{t\geq 0}\overline{\{\bar{\theta}(s): s\geq t\}}$ 
is almost everywhere an internally chain transitive set 
of the differential inclusion
\begin{equation}
\label{slower_di}
\dot{\theta}(t) \in \hat{h}(\theta(t)), 
\end{equation}
where $\hat{h}(\theta)=\{\tilde{h}(\theta,\lambda(\theta),\nu) :  \nu \in D^{(1)}(\theta, \lambda(\theta))\}$. 
We call (\ref{fast}) and (\ref{slower_di}) as the faster and slower d.i. 
to correspond with faster and slower recursions, respectively.
\begin{cor}
\label{main_col}
Under the additional assumption that the inclusion 
\begin{equation}
\dot{\theta}(t)\in \hat{h}(\theta(t))), \nonumber 
\end{equation}
has a global attractor set $A_1$,  
\begin{equation}
(\theta_n, w_n) \to \cup_{\theta^* \in A_1}(\theta^*, \lambda(\theta^*)) \mbox{a.s. as $n \to \infty$.}\nonumber 
\end{equation} 
\end{cor}
\begin{rmk}
In case where the set $D^{(2)}(\theta,w)$ is singleton, we can relax \textbf{(A6)} to local attractors also. The relaxed 
assumption will be
\begin{enumerate}[label=\textbf{(A\arabic*)'}]
\setcounter{enumi}{5} 
\item The function $\hat{g}(\theta, w) = \int g(\theta, w, z)\Gamma^{(2)}_{\theta,w}(dz)$ is Lipschitz 
continuous where $\Gamma^{(2)}_{\theta,w}$ is the only element of $D^{(2)}(\theta,w)$. 
Further, for all $\theta \in \mathbb{R}^d$, the o.d.e
\begin{equation}
\label{cpledode}
\dot{w}(t) = \hat{g}(\theta, w(t))  
\end{equation}
has an asymptotically stable equilibrium $\lambda(\theta)$ with domain of attraction $G_\theta$ 
where $\lambda :  \mathbb{R}^d \to \mathbb{R}^k$ is a Lipschitz map with constant $K$.
Also, assume that $\bigcap_{\theta} G_\theta$ is non-empty.
Moreover, the function $V': G \to 
[0,\infty)$ defined by $V'(\theta,w) = V_\theta(w)$ is continuously differentiable where $V_\theta(.)$ is the Lyapunov function 
(for definition see \cite[Chapter 11.2.3]{borkar1})
for the o.d.e. (\ref{cpledode}) with $\lambda(\theta)$ 
as its attractor, 
and $G=\bigcup_{\theta \in \mathbb{R}^d} \{\{\theta\} \times  G_\theta$\}. 
This extra condition is needed
so that the set graph($\lambda$):=$\{(\theta,\lambda(\theta)): \theta \in \mathbb{R}^d\}$ becomes an asymptotically stable set of the coupled o.d.e 
\begin{equation}
\dot{w}(t) = \hat{g}(\theta(t),w(t)), \dot{\theta}(t) = 0.\nonumber  
\end{equation} 
\end{enumerate}
Note that \textbf{(A6)'} allows multiple attractors (at least one of them 
have to be a point, others can be sets) for the faster o.d.e 
for every $\theta$.
 
Then the statement of Theorem \ref{thm} will be modified as in the following: 
\begin{thm}[Slower timescale result when $\lambda(\theta)$ is a local attractor]Under assumptions \textbf{(A1)-(A5), (A6)'} and  \textbf{(A7)},  
on the event ``$\{w_n\}$ belongs to a compact subset $B$ (depending 
on the sample point) of $\bigcap_{\theta \in \mathbb R^d} G_\theta$ \textbf{eventually}", 
\label{thm_local}
\begin{equation}
(\theta_n, w_n) \to \cup_{\theta^* \in A_0}(\theta^*, \lambda(\theta^*)) \mbox{a.s. as $n \to \infty$.}\nonumber 
\end{equation}
\end{thm}
The requirement  on $\{w_n\}$ is much stronger than the usual local attractor 
statement for Kushner-Clarke lemma \cite[Section II.C]{metivier} which requires 
the iterates to enter a compact set in the domain for attraction of the local attractor 
\textit{infinitely often} only. The reason for imposing this strong assumption is 
that graph($\lambda$) is not a subset of 
any compact set in $\mathbb{R}^{d+k}$, and hence the usual tracking 
lemma kind of arguments do not go through directly. One has to 
relate the limit set of the coupled iterate $(\theta_n, w_n)$ to graph($\lambda$) (See the proof of Lemma \ref{fast_res2}).

\end{rmk}

We present the proof of our main results
in the next section.
\section{Main Results}
\label{mres}

We first discuss an extension of the single time-scale controlled Markov noise framework of \cite{borkar}
under our assumptions to prove our main results. Note that the results of \cite{borkar} 
assume that the state space of the controlled Markov process is Polish which 
may impose  additional conditions that are hard to verify. In this section, other 
than proving our two time-scale results, we prove many 
of the results in \cite{borkar} (which were only
stated there) assuming the state space to be compact.
   
We begin by describing the intuition behind the proof techniques in \cite{borkar}.  
\\ \indent
The space $C([0, \infty); \mathbb{R}^d)$ of continuous
functions from $[0,\infty)$ to $\mathbb{R}^d$ is topologized with the 
coarsest topology such that the map that takes any $f \in C([0, \infty); \mathbb{R}^d)$
to its restriction to $[0,T]$ when viewed as an element of the space $C([0, T]; \mathbb{R}^d)$, 
is continuous for all $T>0$. In other words, 
$f_n \to f$ in this space iff $f_n|_{[0,T]} \to f|_{[0,T]}$.
The other notations used below are the same as those in \cite{borkar,borkar1}. We present a few for easy reference. 
\\ \indent
Consider the single time-scale stochastic approximation recursion with controlled Markov noise:
\begin{equation}
\label{cont_mar}
x_{n+1} = x_n + a(n)\left[h(x_n, Y_n) + M_{n+1}\right]. 
\end{equation}
Define time instants $t(0)=0, t(n)=\sum_{m=0}^{n-1} a(m), n\geq 1$. 
Let $\bar{x}(t), t\geq 0$ be the continuous, piecewise linear 
trajectory defined by $\bar{x}(t(n))=x_n, n\geq 0$, with 
linear interpolation on each interval $[t(n), t(n+1))$, i.e.,  
\begin{equation}
\bar{x}(t) = x_n + (x_{n+1} - x_n)\frac{t-t(n)}{t(n+1)-t(n)}, t \in [t(n), t(n+1)).\nonumber 
\end{equation}
Now, define $\tilde{h}(x,\nu)=\int h(x,z)\nu(dz,U)$ for $\nu \in P(S \times U)$. 
Let $\mu(t), t\geq 0$ be the random process defined by 
$\mu(t)=\delta_{(Y_n,Z_n)}$ for $t \in [t(n), t(n+1)), n\geq 0$, where $\delta_{(y,a)}$ 
is the Dirac measure corresponding to $(y,a)$. Consider 
the non-autonomous o.d.e. 
\begin{equation}
\label{auto}
\dot{x}(t) = \tilde{h}(x(t), \mu(t)). 
\end{equation}
Let $x^s(t), t\geq s$, denote the solution to (\ref{auto}) with $x^s(s)=\bar{x}(s)$, for $s\geq0$. 
Note that $x^s(t), t\in [s, s+T]$ and $x^s(t), t\geq s$ can be viewed as elements of $C([0, T]; \mathbb{R}^d)$ and 
$C([0, \infty); \mathbb{R}^d)$ respectively. With this abuse of notation, it is easy to see that 
$\{x^s(.)|_{[s, s+T]}, s\geq 0\}$ is a 
pointwise bounded and equicontinuous family of functions in $C([0, T]; \mathbb{R}^d)~\forall T >0$.
By Arzela-Ascoli theorem, it is relatively compact. From Lemma~2.2 of \cite{borkar} 
one can see that for all $s(n)\uparrow \infty, \{\bar{x}(s(n)+.)|_{[s(n), s(n)+T]}, \linebreak n\geq 1\}$ 
has a limit point in $C([0, T]; \mathbb{R}^d)~\forall T >0$.
With the above topology for $C([0, \infty); \mathbb{R}^d)$, 
$\{x^s(.), s\geq 0\}$ 
is also relatively compact in $C([0, \infty); \mathbb{R}^d)$ and 
for all $s(n)\uparrow \infty, \{\bar{x}(s(n)+.), n\geq 1\}$ has a limit point in $C([0, \infty); \mathbb{R}^d)$.  
\\ \indent
One can write from (\ref{cont_mar})
the following:
\begin{equation}
\bar{x}(u(n)+t) = \bar{x}(u(n)) + \int_{0}^{t}h(\bar{x}(u(n)+\tau), \nu(u(n)+\tau))d\tau + W^n(t),\nonumber 
\end{equation}
where $u(n)\uparrow \infty, \bar{x}(u(n)+.) \to \tilde{x}(\cdot), \nu(t) = (Y_n,Z_n)$ for $t \in [t(n), t(n+1)), n\geq 0$ 
and $W^n(t) = W(t+u(n)) - W(u(n)), W(t) = W_n + (W_{n+1} - W_n)\frac{t-t(n)}{t(n+1)- t(n)},  
 W_n=\sum_{k=0}^{n-1}a(k)M_{k+1}, 
n\geq 0$. 
From here one cannot directly take limit on both sides as finding limit points of $\nu(s+.)$ as $s \to \infty$
is not meaningful. Now, $h(x,y)=\int h(x,z)\delta_{(y,a)}(dz \times U)$. 
Hence by defining 
$\tilde{h}(x,\rho)=\int h(x,z)\rho(dz)$ and $\mu(t) = \delta_{\nu(t)}$ one can write the above as
\begin{equation}
\label{mu} 
\bar{x}(u(n)+t) = \bar{x}(u(n)) + \int_{0}^{t}\tilde{h}(\bar{x}(u(n)+\tau), \mu(u(n)+\tau))d\tau + W^n(t).
\end{equation}The advantage is that the space $\mathcal{U}$ of measurable functions 
from $[0, \infty)$ to $\mathcal{P}(S \times U)$ is compact metrizable, so sub-sequential limits exist. Note that $\mu(\cdot)$ is not a member
of $\mathcal{U}$, rather we need to fix a sample point, i.e., $\mu(.,\omega) \in \mathcal{U}$. 
For ease of understanding, we abuse the terminology and talk about the limit points $\tilde{\mu}(\cdot)$ of $\mu(s+.)$. 
\\ \indent 
From (\ref{mu}) 
one can infer that the limit $\tilde{x}(\cdot)$ of $\bar{x}(u(n)+.)$ satisfies the o.d.e. 
$\dot{x}(t) = \tilde{h}(x(t), \mu(t))$ with $\mu(\cdot)$ replaced by $\tilde{\mu}(\cdot)$. 
Here each $\tilde{\mu}(t), t \in \mathbb{R}$ in $\tilde{\mu}(\cdot)$ is generated through different limiting
processes each one associated with the compact metrizable space $U_t$ = space of measurable 
functions from $[0,t]$ to $\mathcal{P}(S \times U)$. This will be problematic if we want to further explore the process $\tilde{\mu}(\cdot)$
and convert the non-autonomous o.d.e. into an autonomous one. 
\\ \indent 
Hence the main result is proved using an auxiliary lemma \cite[Lemma~2.3]{borkar} other
than the tracking lemma (Lemma~2.2 of \cite{borkar}). Let $u(n(k)) \uparrow \infty$ be such that
$\bar{x}(u(n(k))+.) \to \tilde{x}(\cdot)$ and 
$\mu(u(n(k))+.) \to \tilde{\mu}(\cdot)$, then using Lemma~2.2 of \cite{borkar} 
one can show that $x^{u(n(k))}(\cdot) \to \tilde{x}(\cdot)$. 
Then the auxiliary lemma shows that the 
o.d.e. trajectory $x^{u(n(k))}(\cdot)$ associated with $\mu(u(n(k))+.)$ 
tracks (in the limit) the o.d.e. trajectory associated with $\tilde{\mu}(\cdot)$. Hence Lemma~2.3 of \cite{borkar} links 
the two limiting processes $\tilde{x}(\cdot)$ and $\tilde{\mu}(\cdot)$ 
in some sense. Note that Lemma~2.3 of \cite{borkar} involves only 
the o.d.e. trajectories, not the interpolated trajectory of the algorithm. 

Consider the iteration 
\begin{equation}
\label{eps2}
\theta_{n+1} = \theta_n + a(n)\left[h(\theta_n, Y_n) + \epsilon_n + M_{n+1}\right],
\end{equation}
where $\epsilon_n \to 0$ and the rest of the notations are same as  \cite{borkar}. 
Specifically, $\{Y_n\}$ is the controlled Markov process driven by $\{\theta_n\}$
and $M_{n+1}, n\geq 0$ is a martingale difference sequence. 
Let $\bar{\theta}(t), t\geq 0$ be the continuous, piecewise linear 
trajectory of (\ref{eps2}) defined by $\bar{\theta}(t(n))=\theta_n, n\geq 0$, with 
linear interpolation on each interval $[t(n), t(n+1))$. Also, 
let $\theta^s(t), t\geq s$, denote the solution to (\ref{auto}) with $\theta^s(s)=\bar{\theta}(s)$, for $s\geq0$. 
\\ \indent
The convergence analysis of (\ref{eps2}) requires some changes  
in Lemma~2.2 and 3.1 of  \cite{borkar}. The modified versions of them are precisely 
the following two lemmas.   
\begin{lem}
\label{track_fast}
For any $T >0$, $\sup_{t\in [s,s+T]}\|\bar{\theta}(t) - \theta^s(t)\| \to 0,$ a.s. as $s\to \infty$. 
\end{lem}
\begin{pf}
The proof follows from the Lemma 2.2 and the remark 3 thereof (p. 144)  of \cite{borkar}.
\end{pf}

Now, $\mu$ can be viewed as a random variable taking values in $\mathcal{U}$ = the space of measurable functions from $[0,\infty)$ to
$\mathcal{P}(S \times U)$. This space is topologized with the coarsest topology such that the map
\begin{equation}
\nu(\cdot) \in \mathcal{U} \to \int_{0}^{T} g(t) \int fd\nu(t)dt \in \mathbb{R}\nonumber 
\end{equation}
is continuous for all $f \in C(S), T>0, g \in L_2[0,T]$. Note that $\mathcal{U}$ is compact metrizable.
\begin{lem}
\label{eps}
Almost surely every limit point of $(\mu(s+.), \bar{\theta}(s+.))$ as $s\to \infty$ is 
of the form $(\tilde{\mu}(\cdot), \tilde{\theta}(\cdot))$
where $\tilde{\mu}(\cdot)$ satisfies $\tilde{\mu}(t) \in D(\tilde{\theta}(t))$ a.e. $t$.  
\end{lem}
\begin{pf}
Suppose that $u(n)\uparrow \infty$, $\mu(u(n)+.) \to \tilde{\mu}(\cdot)$ and $\bar{\theta}(u(n)+.) \to \tilde{\theta}(\cdot)$.
Let $\{f_i\}$ be countable dense in the unit ball of $C(S)$, hence a separating 
class, i.e., $\forall i, \int f_i d\mu = \int f_i d\nu$ 
implies $\mu=\nu$. For 
each $i$,
\begin{equation}
\zeta^i_n = \sum_{m=1}^{n-1}a(m)(f_i(Y_{m+1}) - \int f_i(y)p(dy|Y_m,Z_m, \theta_m)), n \geq 1, \nonumber 
\end{equation}
is a zero-mean martingale with $\mathcal{F}_n = \sigma(\theta_m, Y_m, Z_m, m\leq n)$. 
Moreover, it is a square integrable martingale due to the 
fact that $f_i$'s are bounded and each $\zeta^i_n$ is a finite sum. Its quadratic 
variation process 
\begin{equation}
A_{n}=\sum_{m=0}^{n-1}a(m)^2E[(f_i(Y_{m+1}) - 
\int f_i(y)p(dy|Y_m,Z_m, \theta_m))^2|\mathcal{F}_m] + E[(\zeta^i_0)^2]\nonumber 
\end{equation}
is almost surely convergent. 
By the martingale convergence theorem, 
 $\zeta^i_n, n\geq 0$ converges a.s. for all $i$. As before let $\tau(n,t)=\min\{m \geq n:  t(m) \geq t(n)+t\}$ for $t\geq0, n\geq0$. 
Then as $n\to\infty$,
\begin{equation}
\sum_{m=n}^{\tau(n,t)} a(m)(f_i(Y_{m+1})-\int f_i(y)p(dy|Y_m,Z_m,\theta_m))\to 0,\mbox{ a.s.}\nonumber 
\end{equation}
for $t >0$. 
By our choice of $\{f_i\}$ and the fact that 
$\{a(n)\}$ is an eventually non-increasing sequence (the latter property is used only here and in Lemma \ref{slowmu}), we have
\begin{equation}
\sum_{m=n}^{\tau(n,t)}(a(m) - a(m+1))f_i(Y_{m+1}) \to 0,\mbox{ a.s.}\nonumber 
\end{equation}
From the foregoing, 
\begin{equation}
\sum_{m=n}^{\tau(n,t)} (a(m+1)f_i(Y_{m+1})-a(m)\int f_i(y)p(dy|Y_m,Z_m,\theta_m))\to 0,\mbox{ a.s.}\nonumber 
\end{equation}
for all $t >0$,
which implies
\begin{equation}
\sum_{m=n}^{\tau(n,t)} a(m)(f_i(Y_{m})-\int f_i(y)p(dy|Y_m,Z_m,\theta_m))\to 0,\mbox{ a.s.}\nonumber 
\end{equation}
for all $t >0$ due to the fact that 
$a(n) \to 0$ and $f_i(.)$ are bounded. 
This implies
\begin{equation}
\int_{t(n)}^{t(n)+t}(\int(f_i(z) - \int f_i(y)p(dy|z,a,\hat{\theta}(s)))\mu(s,dzda))ds \to 0,\mbox{ a.s.}\nonumber 
\end{equation}
and that in turn implies
\begin{equation}
\int_{u(n)}^{u(n)+t}(\int(f_i(z) - \int f_i(y)p(dy|z,a,\hat{\theta}(s)))\mu(s,dzda))ds \to 0,\mbox{ a.s.}\nonumber 
\end{equation}
(this is true because $a(n)\to 0$ and 
$f_i(\cdot)$ is bounded)
where $\hat{\theta}(s) = \theta_n$ when $s \in [t(n), t(n+1))$ for $n\geq 0$. Now, one can claim from the above that
\begin{equation}
\int_{u(n)}^{u(n)+t}(\int(f_i(z) - \int f_i(y)p(dy|z,a,\bar{\theta}(s)))\mu(s,dzda))ds \to 0,\mbox{ a.s.}\nonumber 
\end{equation}
This is due to the fact that the map $S \times U \times \mathbb{R}^{d} \ni (z,a,\theta) \to \int f(y)p(dy|z,a,\theta)$ is continuous and hence 
uniformly continuous on the compact set $A = S \times U \times M$ 
where $M$ is the compact set s.t. $\theta_n\in M$ for all $n$.    
Here we also use the fact that $\|\bar{\theta}(s) - \theta_m\|=\|h(\theta_m, Y_m) + \epsilon_m + M_{m+1}\|(s-s_m) \to 0, s\in [t_m, t_{m+1})$ 
as the first two terms inside the norm in the R.H.S are bounded.  
The above convergence is equivalent to 
\begin{equation}
\int_{0}^{t}(\int(f_i(z) - \int f_i(y)p(dy|z,a,\bar{\theta}(s+u(n)))\mu(s+u(n),dzda))ds \to 0,\mbox{ a.s.}\nonumber 
\end{equation} 
Fix a sample point in the probability one set on which the convergence 
above holds for all $i$. Then the convergence 
above leads to 
\begin{equation}
\label{conjun1}
\int_{0}^{t}(\int f_i(z) - \int f_i(y)p(dy|z,a, \tilde{\theta}(s)))\tilde{\mu}(s, dzda)ds =0~\forall i. 
\end{equation}
Here we use one part of the proof from Lemma~2.3 of \cite{borkar} that if $\mu^n(\cdot) \to \mu^{\infty}(\cdot) \in \mathcal{U}$ then 
for any $t>0$,
\begin{equation}
\int_{0}^{t} \int \tilde{f}(s,z,a)\mu^n(s,dzda)ds - \int_{0}^{t} \int \tilde{f}(s,z,a)\mu^{\infty}(s,dzda)ds \to 0,\nonumber 
\end{equation}
for all $\tilde{f} \in 
C([0,t] \times S \times A)$ and the fact that $\tilde{f}_n(s,z,a) = \int f_i(y)p(dy|z,a,\bar{\theta}(s+u(n)))$ 
converges uniformly to $\tilde{f}(s,z,a) = \int f_i(y)p(dy|z,a,\tilde{\theta}(s))$. 
To prove the latter, define $g:C([0,t]) \times [0,t] \times S \times A \to \mathbb{R}$ by $g(\theta(\cdot), s,z,a) 
 = \int f_i(y)p(dy|z,a, \theta(s)))$. To see that $g$ is continuous we need to check that if $\theta_n(\cdot) \to \theta(\cdot)$ 
uniformly and $s(n) \to s$, then $\theta_n(s(n)) \to \theta(s)$. This is because $\|\theta_n(s(n)) - \theta(s)\| = \|\theta_n(s(n)) - \theta(s(n)) + \theta(s(n)) - \theta(s)\| \leq \|\theta_n(s(n)) - \theta(s(n))\| + 
\|\theta(s(n)) - \theta(s)\|$. The first and second terms go to zero 
due to the uniform convergence of $\theta_n(\cdot), n\geq 0$ and continuity of $\theta(\cdot)$ respectively.  
 Let $A = \{\bar{\theta}(u(n)+.)|_{[u(n),u(n)+t]}, n\geq 1\} \cup \tilde{\theta}(\cdot)|_{[0,t]}$. 
$A$ is compact as it is the
union of a sequence of functions and their limit. 
So, $g|_{(A \times [0,t]\times S \times U)}$ is uniformly 
continuous. Then using the same arguments as in Lemma~\ref{upsem} 
we can show equicontinuity of $\{\tilde{f}_n(.,.)\}$, that results in
uniform convergence and thereby (\ref{conjun1}).  
An application of Lebesgue's theorem in conjunction 
with (\ref{conjun1}) shows that
\begin{equation}
\int (f_i(z) - \int f_i(y)p(dy|z,a,\tilde{\theta}(t)))\tilde{\mu}(t, dzda) = 0~\forall i\nonumber 
\end{equation}
for a.e. $t$. 
By our choice of $\{f_i\}$, this leads to 
\begin{equation}
\tilde{\mu}(t, dy \times U) = \int p(dy|z,a,\tilde{\theta}(t))\tilde{\mu}(t, dzda)\nonumber 
\end{equation}
a.e. $t$. Therefore the conclusion follows by disintegrating such measure as the product of marginal on $S$ and
the regular conditional law on $U$ (\cite[p~140]{borkar}). 
\end{pf}

\begin{rmk}
Note that the above invariant distribution does not come ``naturally''; rather it arises from the assumption made to match the natural 
timescale intuition for the controlled Markov noise component, i.e., the slower iterate should see the average effect of the Markov component.
\end{rmk}

The proof of the following lemma, in this case, will be unchanged from its original version, 
so we just mention it for completeness and refer the reader to Lemma 2.3 of \cite{borkar} for its proof.
\begin{lem}
\label{eps3}
Let $\mu^n(\cdot) \to \mu^{\infty}(\cdot) \in \mathcal{U}$. Let $\theta^n(\cdot), n=1, 2, \dots, \infty$ denote 
solutions to (\ref{auto}) corresponding to the case where $\mu(\cdot)$ is replaced by $\mu^n(\cdot)$, for $n=1,2,\dots \infty$. 
Suppose $\theta^n(0) \to \theta^{\infty}(0)$. Then 
\begin{equation}
\lim_{n \to \infty} \sup_{t\in [0, T]}\|\theta^n(t) - \theta^{\infty}(t)\| = 0 \nonumber 
\end{equation}
for every $T >0$. 
\end{lem}
\begin{lem}
Almost surely, $\{\theta_n\}$ converges to an internally 
chain transitive set of the differential inclusion
\begin{equation}
\label{inc}
\dot{\theta}(t) \in \hat{h}(\theta(t)), 
\end{equation}
where $\hat{h}(\theta)=\{\tilde{h}(\theta,\nu) :  \nu \in D(\theta)\}$. 
\end{lem}
\begin{pf}
Lemma~\ref{eps3} shows that every limit point $(\tilde{\mu}(\cdot), \tilde{\theta}(\cdot))$ 
of $(\mu(s+.),\bar{\theta}(s+.))$ as $s\to \infty$ is such that
 $\tilde{\theta}(\cdot)$ satisfies (\ref{auto}) with $\mu(\cdot) = \tilde{\mu}(\cdot)$. 
Hence, $\tilde{\theta}(\cdot)$ is absolutely continuous. Moreover, 
using Lemma~\ref{eps}, one can see that it satisfies (\ref{inc}) a.e. $t$, hence is a solution to the differential inclusion (\ref{inc}).  
Hence the proof follows.
\end{pf}
  
\begin{lem}[Faster timescale result]
 $(\theta_n, w_n) \to \{(\theta, \lambda(\theta)) :  \theta \in \mathbb{R}^d\}$ a.s.
\end{lem}

\begin{pf}
We first rewrite (\ref{eqn1}) as
\begin{equation}
\theta_{n+1} = \theta_n + b(n)\left[\epsilon_n + M^{(3)}_{n+1}\right],\nonumber 
\end{equation}
where $\epsilon_n = \frac{a(n)}{b(n)}h(\theta_n, w_n, Z^{(1)}_n)\to 0$ as $n\to \infty$ a.s. and $M^{(3)}_{n+1} = \frac{a(n)}{b(n)} M^{(1)}_{n+1}$ for $n\geq 0$. 
Let $\alpha_n=(\theta_n, w_n), \alpha=(\theta,w) \in \mathbb{R}^{d+k}, G(\alpha,z)=(0, g(\alpha,z)), \epsilon'_n=(\epsilon_n, 0), 
M^{(4)}_{n+1}= (M^{(3)}_{n+1}, M^{(2)}_{n+1})$. Then one can write 
 (\ref{eqn1}) and (\ref{eqn2}) in the framework of (\ref{eps2}) as
\begin{equation}
\label{stability}
\alpha_{n+1} = \alpha_n + b(n)\left[G(\alpha_n,Z^{(2)}_n) + \epsilon'_n +  M^{(4)}_{n+1}\right],
\end{equation}
with $\epsilon'_n \to 0$ as $n \to \infty$.  
$\alpha_n, n\geq 0$ converges almost surely to an internally chain transitive set of the differential inclusion
\begin{equation}
\dot{\alpha}(t) \in \hat{G}(\alpha(t)),\nonumber 
\end{equation}
where $\hat{G}(\alpha) = \{\tilde{G}(\alpha, \nu) :  \nu \in D^{(2)}(\theta,w)\}$ with $\tilde{G}(\alpha,\nu)=(0,\tilde{g}(\theta,w,\nu))$. 
In other words, 
$(\theta_n, w_n), n\geq 0$ converges to an internally chain transitive set of the differential inclusion
\begin{equation}
\label{coupled_di}
\dot{w}(t) \in \hat{g}(\theta(t), w(t)), \dot{\theta}(t) = 0.\nonumber 
\end{equation}The rest follows from the second part of \textbf{(A6)}.
\end{pf}

\begin{rmk} Under the conditions mentioned in Remark 4 the above 
 faster timescale result should be modified as follows:

\begin{lem}[Faster timescale result when $\lambda(\theta)$ is a local attractor]
\label{fast_res2}
Under assumptions \textbf{(A1) - (A5), (A6)'} and  \textbf{(A7)},  
on the event ``$\{w_n\}$ belongs to a compact subset $B$ (depending 
on the sample point) of $\bigcap_{\theta \in \mathbb R^d} G_\theta$ \textit{eventually}'',
\begin{equation}
(\theta_n, w_n) \to \{(\theta, \lambda(\theta)) :  \theta \in \mathbb{R}^d\} \mbox{~~a.s.} \nonumber 
\end{equation}
\end{lem}
\begin{pf}
Fix a sample point $\omega$.
The proof follows from these observations: 
\begin{enumerate}
 \item continuity of flow for the coupled o.d.e around the initial point,
 \item $\sup_n \|\theta_n\| = M_1 < \infty$,
 \item the fact that the 
set graph($\lambda$) is Lyapunov stable ($V'(.)$ as mentioned in 
\textbf{(A6)'} will be a Lyapunov function for this set), and
  \item the fact that $\bigcap_{t\geq 0} \overline{\bar{\alpha}(s): s \geq t}$ is an internally chain transitive set of 
the coupled o.d.e 
\begin{equation}
\label{copuled_ode}
\dot{w}(t) = \hat{g}(\theta(t),w(t)), \dot{\theta}(t) = 0,  
\end{equation}
where $\bar{\alpha}(.)$ is the interpolated trajectory of the coupled iterate $\{\alpha_n\}$.
\end{enumerate}
As $\{\theta: \|\theta\| \leq M_1\} \times B \subset \bigcup_{\theta \in \mathbb{R}^d} \{\{\theta\} \times G_\theta \}$,  
the first three observations show that for all $\epsilon>0$, there exists a $T_\epsilon >0$ 
such that any o.d.e trajectory for (\ref{copuled_ode}) with starting point 
on the compact set $\{\theta: \|\theta\| \leq M_1\} \times B$ reaches the 
$\epsilon$-neighbourhood of graph($\lambda$) after time $T_\epsilon$.
Further, 
\begin{equation}
\bigcap_{t\geq 0} \overline{\bar{\alpha}(s): s \geq t} \subset  \{\theta: \|\theta\| \leq M_1\} \times B. \nonumber 
\end{equation}
Then one can use 
the last observation by choosing $T > T_{\epsilon}$ to 
show the required convergence to the set  graph($\lambda$).
\end{pf}
\end{rmk}

\begin{rmk}
One interesting question in this context is to analyze whether one can extend the single timescale local attractor 
convergence statements to the two timescale setting under some \textit{verifiable conditions}. More specifically, 
if there is a
global attractor $A_1$ for
\begin{equation}
\dot{\theta}(t) \in \hat{h}(\theta(t)), \nonumber
\end{equation}
then can one provide verifiable conditions to show 
\begin{equation}
P [(\theta_n , w_n ) \to \cup_{\theta \in A_1} (\theta , \lambda(\theta))] > 0. \nonumber
\end{equation}
Here $\lambda(\theta)$ is a local attractor as mentioned in \textbf{(A6)'}.

There are two ways in which this could possibly be tried: 

\begin{enumerate}
 \item Use Theorem \ref{thm_local} where we show that 
 on the event $\{w_n\}$ belongs to a compact subset $B$ (depending 
on the sample point) of $\bigcap_{\theta \in \mathbb R^d} G_\theta$ ``eventually'',
\begin{equation}
(\theta_n, w_n) \to \cup_{\theta^* \in A_1}(\theta^*, \lambda(\theta^*)) \mbox{a.s. as $n \to \infty$,}\nonumber 
\end{equation}which  is an extension of Kushner-Clarke Lemma to the two timescale case.  
Therefore the task would be to impose verifiable assumptions so that
$P$($\{w_n\}$ belongs to a compact subset $B$ (depending 
on the sample point) of $\bigcap_{\theta \in \mathbb R^d} G_\theta$ ``eventually'') $>$ 0. 
In a stochastic approximation scenario it is not immediately clear how one could possibly 
impose verifiable assumptions so that 
such a probabilistic statement becomes true.
\item The second approach would be to extend the analysis of \cite{benaimode, benaim} to the two timescale case.
In our opinion this is very hard as this analysis is based on the attractor introduced by Benaim et al. whereas 
the coupled o.d.e (\ref{copuled_ode}) which tracks the coupled iterate 
$(\theta_n,w_n)$ (therefore the interpolated trajectory of the coupled iterate will be 
an asymptotic pseudo-trajectory \cite{benaimode} for (\ref{copuled_ode}))  has no attractor. The reason is that one cannot obtain a fundamental neighbourhood for sets like 
$\cup_{\theta \in A_1} (\theta, \lambda(\theta))$ as the $\theta$ component will remain constant for any trajectory 
of the above coupled o.d.e.
\end{enumerate}

Thus it is immediately not clear as to how this question can be addressed and this will be an interesting future 
direction.
\end{rmk}

From the faster timescale results we get, $\|w_n - \lambda(\theta_n)\| \to 0$ a.s., i.e, $\{w_n\}$ asymptotically tracks $\{\lambda(\theta_n)\}$ a.s.
\\ \indent 
Now, consider the non-autonomous o.d.e.
\begin{equation}
\label{slow}
\dot{\theta}(t) = \tilde{h}(\theta(t),\lambda(\theta(t)),\mu(t)), 
\end{equation}
where $\mu(t) = \delta_{Z^{(1)}_n,A^{(1)}_n}$ when $t \in [t(n), t(n+1))$ for $n\geq 0$ and $\tilde{h}(\theta,w,\nu)=\int h(\theta,w,z) \nu(dz)$. 
Let $\theta^s(t), t\geq s$ denote the solution to (\ref{slow}) with 
$\theta^s(s) = \bar{\theta}(s)$, for $s \geq 0$. Then
\begin{lem}
\label{track_slow}
For any $T >0, \sup_{t\in [s,s+T]}\|\bar{\theta}(t) - \theta^s(t)\| \to 0,$ a.s.  
\end{lem}
\begin{pf}
The slower recursion corresponds to
\begin{equation}
\theta_{n+1} = \theta_n + a(n)\left[h(\theta_n, w_n, Z^{(1)}_n) + M^{(1)}_{n+1}\right].\nonumber  
\end{equation}
Let $t(n+m) \in [t(n), t(n) + T]$.  Let $[t] =\max\{t(k) :  t(k) \leq t\}$. Then by construction,
\begin{align}
\bar{\theta}(t(n+m)) &= \bar{\theta}(t(n)) + \sum_{k=0}^{m-1} a(n+k)h(\bar{\theta}(t(n+k)), w_{n+k}, Z^{(1)}_{n+k}) + \delta_{n, n+m} \nonumber\\
                &= \bar{\theta}(t(n)) + \sum_{k=0}^{m-1} a(n+k)h(\bar{\theta}(t(n+k)), \lambda(\bar{\theta}(t(n+k))), Z^{(1)}_{n+k})\nonumber\\
                &+\sum_{k=0}^{m-1} a(n+k)(h(\bar{\theta}(t(n+k)), w_{n+k}, Z^{(1)}_{n+k})- h(\bar{\theta}(t(n+k)), \lambda(\theta_{n+k}), Z^{(1)}_{n+k}))\nonumber\\
                &+ \delta_{n, n+m},\nonumber
\end{align}
where $\delta_{n, n+m}=\zeta_{n+m}- \zeta_{n}$ with  $\zeta_{n} = \sum_{m=0}^{n-1}a(m)M^{(1)}_{m+1}, n \geq 1$.
\begin{align}
 \theta^{t(n)}(t(m+n)) &= \bar{\theta}(t(n)) + \int_{t(n)}^{t(n+m)} \tilde{h}(\theta^{t(n)}(t), \lambda(\theta^{t(n)}(t)),  \mu(t))dt\nonumber\\
                  &= \bar{\theta}(t(n)) + \sum_{k=0}^{m-1} a(n+k)h(\theta^{t(n)}(t(n+k)), \lambda(\theta^{t(n)}(t(n+k))), Z^{(1)}_{n+k})\nonumber\\
                  &+ \int_{t(n)}^{t(n+m)} (h(\theta^{t(n)}(t), \lambda(\theta^{t(n)}(t), \mu(t))) - h(\theta^{t(n)}([t]), \lambda(\theta^{t(n)}([t]), \mu([t]))))dt.\nonumber
\end{align}
Let $t(n) \leq t \leq t(n+m)$. Now, if $0 \leq k \leq (m-1)$ and $t \in (t(n+k), t(n+k+1)],$
\begin{align}
\|\theta^{t(n)}(t)\| &\leq \|\bar{\theta}(t(n)\| + \|\int_{t(n)}^{t} \tilde{h}(\theta^{t(n)}(\tau), \lambda(\theta^{t(n)}(\tau)), \mu(\tau))d\tau\|\nonumber\\
                &\leq \|\theta_n\| + \sum_{l=0}^{k-1} \int_{t(n+l)}^{t(n+l+1)} (\|h(0,0,Z^{(1)}_{n+l})\|+ L^{(1)}(\|\lambda(0)\|+(K+1)\|\theta^{t(n)}(\tau)\|))d\tau\nonumber\\
                & +\int_{t(n+k)}^t(\|h(0,0,Z^{(1)}_{n+k})\|+ L^{(1)}(\|\lambda(0)\|+(K+1)\|\theta^{t(n)}(\tau)\|))d\tau\nonumber\\
                &\leq C_0 + (M+L^{(1)}\|\lambda(0)\|)T + L^{(1)}(K+1)\int_{t(n)}^{t}\|\theta^{t(n)}(\tau)\|d\tau,\nonumber 
\end{align}
where $C_0 = \sup_n \|\theta_n\| < \infty, \sup_{z \in S^{(1)}}\|h(0,0,z)\| = M$. 
By Gronwall's inequality, it follows that 
\begin{equation}
\|\theta^{t(n)}(t)\| \leq (C_0 + (M+L^{(1)}\|\lambda(0)\|)T)e^{L^{(1)}(K+1)T}.\nonumber 
\end{equation}
\begin{align}
\|\theta^{t(n)}(t) - \theta^{t(n)}(t(n+k))\| &\leq \int_{t(n+k)}^t\|h(\theta^{t(n)}(s), \lambda(\theta^{t(n)}(s)), Z^{(1)}_{n+k})\|ds\nonumber\\
                                   &\leq (\|h(0,0,Z^{(1)}_{n+k})\|+ L^{(1)}\|\lambda(0)\|)(t-t(n+k))\nonumber\\
                                   &+L^{(1)}(K+1)\int_{t(n+k)}^t\|\theta^{t(n)}(s)\|ds\nonumber\\
                                   &\leq C_Ta(n+k),\nonumber
\end{align}
where $C_T=(M+ L^{(1)}\|\lambda(0)\|) + L^{(1)}(K+1)(C_0 + (M+L^{(1)}\|\lambda(0)\|)T)e^{L^{(1)}(K+1)T}.$
Thus, 
\begin{align}
&\|\int_{t(n)}^{t(n+m)} (h(\theta^{t(n)}(t), \lambda(\theta^{t(n)}(t)), \mu(t)) - h(\theta^{t(n)}([t]), \lambda(\theta^{t(n)}([t])), \mu([t])))dt\|\nonumber\\
&\leq \sum_{k=0}^{m-1}\int_{t(n+k)}^{t(n+k+1)}\|h(\theta^{t(n)}(t), \lambda(\theta^{t(n)}(t)), Z^{(1)}_{n+k}) - h(\theta^{t(n)}([t]), \lambda(\theta^{t(n)}([t])), Z^{(1)}_{n+k})\|dt\nonumber\\
&\leq L\sum_{k=0}^{m-1}\int_{t(n+k)}^{t(n+k+1)}\|\theta^{t(n)}(t) - \theta^{t(n)}(t(n+k))\|dt\nonumber\\
&\leq C_TL\sum_{k=0}^{m-1}a(n+k)^2\nonumber\\
&\leq C_TL\sum_{k=0}^{\infty}a(n+k)^2 \to 0 \mbox{~as $n \to \infty$},\nonumber
\mbox{where $L=L^{(1)}(K+1).$}
\end{align}
Hence
\begin{align}
\|\bar{\theta}(t(n+m))-\theta^{t(n)}(t(n+m))&\leq L\sum_{k=0}^{m-1}a(n+k)\|\bar{\theta}(t(n+k)) - \theta^{t(n)}(t(n+k))\|\nonumber\\
& + C_TL\sum_{k=0}^{\infty}a(n+k)^2 + \sup_{k\geq 0}\|\delta_{n,n+k}\|\nonumber\\
& + L^{(1)}\sum_{k=0}^{m-1}a(n+k)\|w_{n+k} -\lambda(\theta_{n+k})\|\nonumber\\
&\leq L\sum_{k=0}^{m-1}a(n+k)\|\bar{\theta}(t(n+k)) - \theta^{t(n)}(t(n+k))\|\nonumber\\
& + C_TL\sum_{k=0}^{\infty}a(n+k)^2 + \sup_{k\geq 0}\|\delta_{n,n+k}\|\nonumber\\
& +  L^{(1)}T\sup_{k\geq 0} \|w_{n+k} - \lambda(\theta_{n+k})\|,\nonumber \mbox{ a.s.}
\end{align}
Define
\begin{equation}
K_{T,n} =  C_TL\sum_{k=0}^{\infty}a(n+k)^2 + \sup_{k\geq 0}\|\delta_{n,n+k}\|
 + L^{(1)}T\sup_{k\geq 0} \|w_{n+k} - \lambda(\theta_{n+k})\|. \nonumber 
\end{equation}
Note that $K_{T,n} \to 0$ a.s. 
The remainder of the proof follows in the exact 
same manner as the tracking lemma, see Lemma 1, Chapter 2 of \cite{borkar1}. 
\end{pf}

\begin{lem}
\label{ode}
Suppose, $\mu^n(\cdot) \to \mu^{\infty}(\cdot) \in U^{(1)}$. Let $\theta^n(\cdot), n=1, 2, \dots, \infty$ 
denote solutions to (\ref{slow}) corresponding to the case where $\mu(\cdot)$ 
is replaced by $\mu^n(\cdot)$, for $n=1,2, \dots, \infty$. 
Suppose $\theta^n(0) \to \theta^{\infty}(0)$. Then 
\begin{equation}
\lim_{n \to \infty} \sup_{t\in [0, T]}\|\theta^n(t) - \theta^{\infty}(t)\| \to 0 \nonumber 
\end{equation}
for every $T >0$. 
\end{lem} 
\begin{pf}
It is shown in  Lemma~2.3 of \cite{borkar} that
\begin{equation}
\int_{0}^{t}\int \tilde{f}(s,z)\mu^{n}(s,dz)ds - \int_{0}^{t}\int \tilde{f}(s,z)\mu^{\infty}(s,dz)ds \to 0 \nonumber 
\end{equation}
for any $\tilde{f} \in C([0,T]\times S)$. 
Using this, one can see that
\begin{equation}
\|\int_{0}^{t} (\tilde{h}(\theta^{\infty}(s),\lambda(\theta^{\infty}(s)),  \mu^n(s)) - 
\tilde{h}(\theta^{\infty}(s), \lambda(\theta^{\infty}(s)), \mu^{\infty}(s)))ds \| \to 0.\nonumber 
\end{equation}
This follows because $\lambda$ is continuous and $h$ is jointly continuous in its arguments.
As a function of $t$, the integral on the left is equicontinuous and pointwise bounded. By the Arzela-Ascoli theorem, this 
convergence must in fact be uniform for $t$ in a compact set. Now for $t>0$, 
\begin{align}
&\|\theta^n(t)-\theta^{\infty}(t)\|\nonumber\\ 
&\leq \|\theta^n(0) - \theta^{\infty}(0)\| + \int_{0}^{t} \|\tilde{h}(\theta^n(s), \lambda(\theta^n(s)),\mu^n(s)) - \tilde{h}(\theta^{\infty}(s), \lambda(\theta^{\infty}(s)), \mu^{\infty}(s))\|ds\nonumber\\
&\leq \|\theta^n(0) - \theta^{\infty}(0)\| + \int_{0}^{t} (\|\tilde{h}(\theta^n(s), \lambda(\theta^n(s)),\mu^n(s)) - \tilde{h}(\theta^{\infty}(s), \lambda(\theta^{\infty}(s)), \mu^{n}(s))\|)ds\nonumber\\ 
&+ \int_{0}^{t} (\|\tilde{h}(\theta^{\infty}(s),\lambda(\theta^{\infty}(s)), \mu^{n}(s)) - \tilde{h}(\theta^{\infty}(s),\lambda(\theta^{\infty}(s)), \mu^{\infty}(s))\|)ds.\nonumber
\end{align}
Now, using the fact that $\lambda$ is Lipschitz with constant $K$ the remaining 
part of the proof follows in the same manner as Lemma~2.3 of \cite{borkar}. 
\end{pf}

Note that Lemma~\ref{ode} shows that every limit point $(\tilde{\mu}(\cdot), \tilde{\theta}(\cdot))$ 
of $(\mu(s+.),\bar{\theta}(s+.))$ as $s\to \infty$ is such that
 $\tilde{\theta}(\cdot)$ satisfies (\ref{slow}) with $\mu(\cdot) = \tilde{\mu}(\cdot)$. 
\begin{lem}
\label{slowmu}
Almost surely every limit point of $(\mu(s+.),\bar{\theta}(s+.))$ as 
$s \to \infty$ is of the form $(\tilde{\mu}(\cdot), \tilde{\theta}(\cdot))$, where $\tilde{\mu}(\cdot)$ 
satisfies $\tilde{\mu}(t) \in D^{(1)}(\tilde{\theta}(t), \lambda(\tilde{\theta}(t)))$. 
\end{lem}
\begin{pf}
Suppose that $u(n)\uparrow \infty$, $\mu(u(n)+.) \to \tilde{\mu}(\cdot)$ and $\bar{\theta}(u(n)+.) \to \tilde{\theta}(\cdot)$.
Let $\{f_i\}$ be countable dense in the unit ball of $C(S)$, hence it is a separating 
class, i.e., for all $i$, $\int f_i d\mu = \int f_i d\nu$ 
implies $\mu=\nu$. 
For 
each $i$,
\begin{equation}
\zeta^i_n = \sum_{m=1}^{n-1}a(m)(f_i(Z^{(1)}_{m+1}) - \int f_i(y)p(dy|Z^{(1)}_m, A^{(1)}_m, \theta_m, w_m)),\nonumber 
\end{equation}
is a zero-mean martingale with $\mathcal{F}_n = \sigma(\theta_m, w_m, Z^{(1)}_m,A^{(1)}_m, m\leq n), n\geq 1$. 
Moreover, it is a square-integrable martingale due to the 
fact that $f_i$'s are bounded and each $\zeta^i_n$ is a finite sum. Its quadratic variation process 
\begin{equation}
A_n=\sum_{m=0}^{n-1}a(m)^2  E[(f_i(Z^{(1)}_{m+1}) - 
\int f_i(y)p(dy|Z^{(1)}_m, A^{(1)}_m, \theta_m, w_m))^2|\mathcal{F}_m] + E[(\zeta^i_0)^2]\nonumber 
\end{equation}
is almost surely convergent. 
By the martingale convergence theorem, $\{\zeta^i_n\}$  
converges a.s. Let $\tau(n,t)=\min\{m \geq n:  
t(m) \geq t(n)+t\}$ for $t\geq0, n\geq0$. Then as $n\to\infty$,
\begin{equation}
\sum_{m=n}^{\tau(n,t)} a(m)(f_i(Z^{(1)}_{m+1})-\int f_i(y)p(dy|Z^{(1)}_m, A^{(1)}_m,\theta_m, w_m))\to 0,\mbox{ a.s.,}\nonumber 
\end{equation}
for $t >0$. 
By our choice of $\{f_i\}$ and the fact that 
$\{a(n)\}$ are eventually non-increasing,
\begin{equation}
\sum_{m=n}^{\tau(n,t)}(a(m) - a(m+1))f_i(Z^{(1)}_{m+1}) \to 0,\mbox{a.s.}\nonumber 
\end{equation} 
Thus,
\begin{equation}
\sum_{m=n}^{\tau(n,t)} a(m)(f_i(Z^{(1)}_m)-\int f_i(y)p(dy|Z^{(1)}_m, A^{(1)}_m,\theta_m,w_m))\to 0,\mbox{ a.s.}\nonumber 
\end{equation}
which implies
\begin{equation}
\int_{t(n)}^{t(n)+t}(\int(f_i(z) - \int f_i(y)p(dy|z,a,\hat{\theta}(s),\hat{w}(s)))\mu(s,dzda))ds \to 0,\mbox{ a.s.}\nonumber 
\end{equation} 
Recall that $u(n)$ can be any general sequence other than $t(n)$. Therefore
\begin{equation}
\int_{u(n)}^{u(n)+t}(\int(f_i(z) - \int f_i(y)p(dy|z,a,\hat{\theta}(s),\hat{w}(s)))\mu(s,dzda))ds \to 0,\mbox{ a.s.,}\nonumber 
\end{equation}
(this follows from the fact that $a(n)\to 0$ and $f_i$'s are bounded)
where $\hat{\theta}(s) = \theta_n$ and $\hat{w}(s) = w_n$ when $s \in [t(n), t(n+1)), n\geq 0$. Now, one can claim from the above that
\begin{equation}
\int_{u(n)}^{u(n)+t}(\int(f_i(z) - \int f_i(y)p(dy|z,a,\bar{\theta}(s), \lambda(\bar{\theta}(s))))\mu(s,dzda))ds \to 0,\mbox{ a.s.}\nonumber 
\end{equation}
This is due to the fact that the map $S^{(1)} \times U^{(1)} \times \mathbb{R}^{d+k} \ni (z,a,\theta,w) \to \int f_i(y)p(dy|z,a,\theta,w)$ is continuous and hence 
uniformly continuous on the compact set $A =  \linebreak S^{(1)} \times U^{(1)} \times M_1 \times M_2$ 
where $M_1$ is the compact set s.t. $\theta_n \in M_1$ for all $n$ and 
 $M_2=\linebreak \{w :  \|w\| \leq \max(\sup\|w_n\|, K')\}$
where $K'$ is the bound for the compact set $\lambda(M_1)$.    
Here we also use the fact that $\|w_m - \lambda(\bar{\theta}(s))\|\to 0$ for $s\in [t_m, t_{m+1})$ as $\lambda$ is Lipschitz and $\|w_m -\lambda(\theta_m)\| \to 0$. 
The above convergence is equivalent to
\begin{equation}
\int_{0}^{t}(\int(f_i(z) - \int f_i(y)p(dy|z,a,\bar{\theta}(s+u(n)), \lambda(\bar{\theta}(s+u(n)))))\mu(s+u(n),dzda))ds \to 0\mbox{ a.s.}\nonumber 
\end{equation} 
Fix a sample point in the probability one set on which the convergence above holds for all $i$. Then the convergence 
above leads to 
\begin{equation}
\label{conjun}
\int_{0}^{t}(\int f_i(z) - \int f_i(y)p(dy|z,a, \tilde{\theta}(s), \lambda(\tilde{\theta}(s))))\tilde{\mu}(s, dzda)ds =0~\forall i. 
\end{equation}
For showing the above, we use one part of the proof from  Lemma~2.3 of \cite{borkar} that if $\mu^n(\cdot) \to \mu^{\infty}(\cdot) \in \mathcal{U}$ then 
for any $t$,
\begin{equation}
\int_{0}^{t} \int \tilde{f}(s,z,a)\mu^n(s,dzda)ds - \int_{0}^{t} \int \tilde{f}(s,z,a)\mu^{\infty}(s,dzda)ds \to 0\nonumber 
\end{equation}
for all $\tilde{f} \in 
C([0,t] \times S^{(1)} \times U^{(1)})$. In addition, we make use of the fact that $\tilde{f}_n(s,z,a) =  \linebreak \int f_i(y)p(dy|z,a,\bar{\theta}(s+u(n)), \lambda(\bar{\theta}(s+u(n))))$ 
converges uniformly to  $\tilde{f}(s,z,a) =  \int f_i(y)p(dy|z,a,\tilde{\theta}(s), \lambda(\tilde{\theta}(s)))$. To prove this, define
$g :C([0,t]) \times [0,t] \times S^{(1)} \times U^{(1)} \to \mathbb{R}$ by  $g(\theta(\cdot), s,z,a) 
 = \int f_i(y)p(dy|z,a, \theta(s),\lambda(\theta(s)))$. 
Let $A' = \{\bar{\theta}(u(n)+.)|_{[u(n),u(n)+t]}, n \geq 1\} \cup \tilde{\theta}(\cdot)|_{[0,t]}$. 
Using the same argument as in Lemma~\ref{eps} and \textbf{(A6)}, i.e., $\lambda$ is Lipschitz (the latter 
helps to claim that if $\theta_n(\cdot) \to \theta(\cdot)$
uniformly then $\lambda(\theta_n(\cdot)) \to \lambda(\theta(\cdot))$ uniformly),  
it can be seen that $g$ is continuous. Then $A'$ is compact  
as it is a union of a sequence of functions and its limit. 
So, $g|_{(A'\times [0,t] \times S^{(1)} \times U^{(1)})}$ is uniformly 
continuous. Then a similar argument as in Lemma~\ref{upsem} shows 
equicontinuity of $\{\tilde{f}_n(.,.)\}$ that results in 
uniform convergence and thereby (\ref{conjun}).
An application of Lebesgue's theorem in conjunction 
with (\ref{conjun}) shows that
\begin{equation}
\int (f_i(z) - \int f_i(y)p(dy|z,a,\tilde{\theta}(t), \lambda(\tilde{\theta}(t)))\tilde{\mu}(t, dzda) = 0~\forall i\nonumber 
\end{equation}
for a.e. $t$. 
By our choice of $\{f_i\}$, this leads to
\begin{equation}
\tilde{\mu}(t, dy \times U^{(1)}) = \int p(dy|z,a,\tilde{\theta}(t), \lambda(\tilde{\theta}(t)))\tilde{\mu}(t, dzda),\nonumber 
\end{equation}
a.e. $t$.     
\end{pf}

Lemma~\ref{ode} shows that every limit point $(\tilde{\mu}(\cdot), \tilde{\theta}(\cdot))$ of $(\mu(s+.),\bar{\theta}(s+.))$ as $s\to \infty$ is such that
$\tilde{\theta}(\cdot)$ satisfies (\ref{slow}) with $\mu(\cdot) = \tilde{\mu}(\cdot)$. Hence, $\tilde{\theta}(\cdot)$ is absolutely continuous. Moreover, 
using Lemma~\ref{slowmu}, one can see that it satisfies (\ref{slower_di}) a.e. $t$, hence is a solution to the differential inclusion (\ref{slower_di}).  
\\ \indent

\begin{pf}[Proof of Theorem \ref{thm} and \ref{thm_local}] 
From the previous three lemmas it is easy to see that 
$A_0 = \cap_{t\geq 0}\overline{\{\bar{\theta}(s): s\geq t\}}$ is almost everywhere an internally chain transitive set of (\ref{slower_di}).
\end{pf}

\begin{pf}[Proof of Corollary \ref{main_col}]
Follows directly from Theorem \ref{thm} and Lemma~\ref{ga}.
\end{pf}

\section{Discussion on the assumptions: Relaxation of (A2)} 
\label{relax}
We discuss relaxation of the uniformity of the Lipschitz constant 
w.r.t state of the controlled Markov process for the 
vector field. The modified 
assumption here is
\begin{enumerate}[label=\textbf{(A\arabic*)'}]
\setcounter{enumi}{1}
 \item $h :  \mathbb{R}^{d+k} \times S^{(1)} \to \mathbb{R}^d$ is  
jointly continuous as well as Lipschitz in its first two arguments with the third argument fixed to same value and 
Lipschitz constant is a function of this value. The latter condition means that
\begin{equation}
\forall z^{(1)} \in S^{(1)}, \|h(\theta, w, z^{(1)}) - h(\theta', w', z^{(1)})\| \leq L^{(1)}(z^{(1)})(\|\theta-\theta'\| + \|w - w'\|).\nonumber
\end{equation}
A similar condition holds for $g$ where the Lipschitz constant is $L^{(2)}: S^{(2)} \to \mathbb{R}^+$.
\end{enumerate}
Note that this allows $L^{(i)}(.)$ to be an unbounded measurable function making it discontinuous due to \textbf{(A1)}.
The straightforward solution for implementing this is to additionally assume the following:
\begin{enumerate}[label=\textbf{(A\arabic*)}]
\setcounter{enumi}{7}
\item $\sup_n L^{(i)}(Z^{(i)}_n) < \infty$ a.s.  
\end{enumerate}
still allowing $L^{(i)}(.)$ to be an unbounded function.
As all our proofs in Section \ref{mres} are shown for every sample point of a probability 1 set, 
our proofs will go through. In the following we give such an example for the case where the Markov 
process is uncontrolled. 

It is enough to consider examples with locally compact $S^{(i)}$ (because then we can take the standard one-point compactification 
and define $L^{(i)}$ arbitrarily at the extra point).

Let $S^{(i)}=\Bbb Z$ and let $Z^{(i)}_n, n \geq 0$ be the Markov Chain on $\Bbb Z$ starting at $0$ with transition probabilities
$p(n,n+1)=p$ and $p(n,n-1)=1-p$. We assume $1/2 < p < 1$. Let $L^{(i)}(n) = \big( \frac{1-p}{p} \big)^n$.

Note that $Z^{(i)}_n,n \geq 0$ is a transient Markov Chain with $Z^{(i)}_n \to +\infty$ a.s. 
From this it follows that $\inf_n Z^{(i)}_n>-\infty$, and 
thus $\sup_n L^{(i)}(Z^{(i)}_n)< \infty$ almost surely. 
It follows that $(L^{(i)}(Z^{(i)}_n))_{n \in \Bbb N}$ is a bounded sequence with probability $1$, but this bound is 
clearly not deterministic since there is a non-zero probability that the sample path reaches large negative values.

However in the following we discuss on the idea of using moment assumptions to
analyze the convergence of single timescale controlled Markov noise framework of \cite{borkar}. 
We show that the iterates (\ref{eps2}) (with $\epsilon_n=0$)  converge to an internally 
chain transitive set of the o.d.e. (\ref{auto}). For this we prove 
Lemma \ref{track_fast} under the following assumptions: 
For all $T >2, i=1,2$,  
\begin{enumerate}[label=\textbf{(S\arabic*)}]
 \item The controlled Markov process $Y_n$ as described in \cite{borkar} takes values in a compact metric space. 
\item For all $n>0$, $0< a(n) \leq 1$, $\sum_n a(n) = \infty$, $\sum_n a(n)^2 < \infty$ and $a(n+1)\leq a(n), n\geq 0$. 
\item $h :  \mathbb{R}^{d} \times S \to \mathbb{R}^d$ Lipschitz in its first argument w.r.t the second. 
The condition means that
\begin{equation}
\forall z \in S, \|h(\theta, z) - h(\theta', z)\| \leq L(z)(\|\theta-\theta'\|).\nonumber
\end{equation}
 \item Let  $\phi(n,T) = \max(m: a(n) + a(n+1) + \dots + a(n+m) \leq T)$ with the bound depending on $T$. 
Then \begin{equation}
        \sup_n E\left[\left(\sup_{0 \leq m \leq \phi(n,T)}L(Y_{n+m})\right)^{16}\right] < \infty. \nonumber
       \end{equation}
 \item \begin{equation}
        \sup_n E\left[e^{8\sum_{m=0} ^ {\phi(n,T)} a(n+m) L(Y_{n+m})}\right] < \infty. \nonumber
       \end{equation}
Note that \textbf{(S4)} and \textbf{(S5)} are trivially satisfied in the case when $L(z)= L$ for all $z \in S$ i.e. the 
case of Section \ref{sec_def}.
\begin{rmk}
As long as one can prove Lemma \ref{track_fast} for all $T >2$ it will hold for all $T>0$, thus one can combine \textbf{(S4)} and 
\textbf{(S5)} into the following assumption:
\begin{equation}
\sup_n E\left[e^{8T\sup_{0 \leq m \leq \phi(n,T)} L(Y_{n+m})}\right] < \infty. \nonumber 
\end{equation} 
As an instance where such an assumption is verified, consider the Markov process of \cite[Eqn. (3.4)]{metivier} defined by 
\begin{equation}
Y_{n+1}= A(\theta_n) Y_n + B(\theta_n) W_{n+1} \nonumber
\end{equation}
where $A(\theta), B(\theta), \theta \in \mathbb{R}^d$, are
$k \times k$-matrices and $(W_n)_{n\geq O}$ are independent and 
identically distributed $\mathbb{R}^k$-valued random variables.  
Assume that the following 
conditions hold true for all $x,y \in S$:
\begin{enumerate}
 \item $L(Y_n)$ is a non-decreasing sequence.
 \item For $r>0, R>0$, \begin{equation}
\sup_{\|\theta\| \leq R} e^{rL(A(\theta) x + B(\theta) y)} \leq L_R {\alpha_R}^r e^{r L(x)} + M_R e^{C_R L(y)} \nonumber 
\end{equation}
for some $C_R, M_R, L_R >0$ and $\alpha_R < 1$.
\end{enumerate}Then  
\begin{align}
\nonumber
&E\left[e^{rL(Y_n)}|Y_{n-1} = x, \theta_{n-1} = \theta\right] \\ \nonumber
&\leq \int e^{rL\left(A(\theta)x + B(\theta) y\right)} \mu_n (dy) \\ \nonumber
&\leq L_R {\alpha_R}^r  e^{rL(x)} + M_R E\left[e^{C_R L(W_n)}\right] \\ \nonumber
&=L_R{\alpha_R}^r  e^{r L(x)} + K_R,\nonumber
\end{align}
with $K_R = M_R E\left[e^{C_R L(W_n)}\right]$ (this follows from the fact that $W_n$ are i.i.d 
if we assume that $E\left[e^{C_R L(W_1)}\right] < \infty$). Choosing large values of $r$, one can show that 
\begin{equation}
E\left[e^{rL(Y_n)}|Y_{n-1} = x, \theta_{n-1} = \theta\right] \leq \beta_R e^{r L(x)}+ K_R \nonumber
\end{equation}
where $\beta_R = L_R{\alpha_R}^r < 1$.
Using the above, for large $r$ 
\begin{align}
E\left[e^{r L(Y_n)}\right] = E\left[E\left[e^{r L(Y_n)}| Y_{n-1}, \theta_{n-1}\right]\right] \leq \beta_R  E\left[e^{rL(Y_{n-1})}\right] + K_R, \nonumber
\end{align}
which shows that 
\begin{equation}
\sup_n E\left[e^{rL(Y_n)}\right] < \infty.\nonumber  
\end{equation} Choosing $r > 8T$,  
\begin{equation}
\sup_n E\left[e^{8TL(Y_n)}\right] < \infty.\nonumber  
\end{equation}
 Note that this is a much weaker assumption that \textbf{(A8)}.
\end{rmk}

 \item The noise sequence $M_n, n \geq 0$ (need not be a martingale difference sequence) 
satisfies the following condition \begin{equation}
        \sup_n E\left[\left(\sum_{m=0} ^ {\phi(n,T)} \|M_{n+m+1}\|\right)^4\right] < \infty. \nonumber
       \end{equation}
\item $\sup_n \|\theta_n\| < \infty$.
\end{enumerate}
With the above assumptions we prove the following tracking lemma:
\begin{lem}
\label{track}
For any $T >0, \sup_{t\in [s,s+T]}\|\bar{\theta}(t) - \theta^s(t)\| \to 0,$ a.s.  
\end{lem}

\begin{pf}
Let $t(n) \leq t \leq t(n+m)$. Now, if $0 \leq k \leq (m-1)$ and $t \in (t(n+k), t(n+k+1)],$
\begin{align}
\|\theta^{t(n)}(t)\| &\leq \|\bar{\theta}(t(n)\| + \|\int_{t(n)}^{t} \tilde{h}(\theta^{t(n)}(\tau), \mu(\tau))d\tau\|\nonumber\\
                &\leq \|\theta_n\| + \sum_{l=0}^{k-1} \int_{t(n+l)}^{t(n+l+1)} (\|h(0,Y_{n+l})\|+ L(Y_{n+l})\|\theta^{t(n)}(\tau)\|))d\tau\nonumber\\
                & +\int_{t(n+k)}^t(\|h(0,Y_{n+k})\|+ L(Y_{n+k})\|\theta^{t(n)}(\tau)\|))d\tau\nonumber\\
                &\leq C_0 + MT+ \int_{t(n)}^{t} L(Y(\tau)) \|\theta^{t(n)}(\tau)\|d\tau \nonumber
\end{align}
where $Y(\tau) = Y_n$ if $\tau \in [t(n),t(n+1))$. Then it follows from an application of Gronwall inequality that
\begin{equation}
\|\theta^{t(n)}(t)\| \leq C e^{\int_{t(n)}^{t} L(Y(\tau)) d\tau} \mbox{~~a.e. $t$}\nonumber
\end{equation}
where $C=C_0 + MT$.
Next, 
\begin{align}
\|\theta^{t(n)}(t) - \theta^{t(n)}(t(n+k))\| &\leq \int_{t(n+k)}^t\|h(\theta^{t(n)}(s), Y_{n+k})\|ds\nonumber\\
                                   &\leq \|h(0,Y_{n+k})\|(t-t(n+k))+ L(Y_{n+k})\int_{t(n+k)}^t\|\theta^{t(n)}(s)\|ds\nonumber\\
                                   &\leq M a(n+k) + C L(Y_{n+k})\int_{t(n+k)}^te^{\int_{t(n)}^s L(Y(\tau)) d\tau} ds.\nonumber
\end{align}
Then 
\begin{align}
&\|\int_{t(n)}^{t(n+m)} (h(\theta^{t(n)}(t), \mu(t)) - h(\theta^{t(n)}([t]), \mu([t])))dt\|\nonumber\\
&\leq \sum_{k=0}^{m-1}\int_{t(n+k)}^{t(n+k+1)}\|h(\theta^{t(n)}(t), Y_{n+k}) - h(\theta^{t(n)}([t]), Y_{n+k})\|dt\nonumber\\
&\leq \sum_{k=0}^{m-1}L(Y_{n+k}) \int_{t(n+k)}^{t(n+k+1)} \|\theta^{t(n)}(t) - \theta^{t(n)}(t(n+k))\|dt\nonumber\\
&\leq \sum_{k=0}^{m-1} c_k\nonumber
\end{align}
where 
\begin{equation}
c_k = L(Y_{n+k})a(n+k)^2\left[M + CL(Y_{n+k}) e^{\sum_{i=0}^{k} a(n+i) L(Y_{n+i})}\right].\nonumber 
\end{equation}
\begin{align}
\|\bar{\theta}(t(n+m))-\theta^{t(n)}(t(n+m))\| &\leq \sum_{k=0}^{m-1}L(Y_{n+k})a(n+k)\|\bar{\theta}(t(n+k)) - \theta^{t(n)}(t(n+k))\|\nonumber\\
& + \sum_{k=0}^{m-1} c_k + \|\delta_{n,n+m}\|,\nonumber 
\end{align}
where $\delta_{n,n+m}=\sum_{k=n}^{n+m-1}a(k)M_{k+1}$.

Therefore using discrete Gronwall inequality we get
\begin{equation}
\|\bar{\theta}(t(n+m))-\theta^{t(n)}(t(n+m))\| \leq r(m,n)  e^{\sum_{k=0}^{m-1} a(n+k) L(Y_{n+k})}\nonumber
\end{equation}
where $r(m,n) = \sum_{k=0}^{m-1} (c_k + a(n+k) \|M_{n+k+1}\|)$.

Now, for some $\lambda \in [0,1]$, 
\begin{align}
&\|\theta^{t(n)}(t) - \bar{\theta}(t)\| \nonumber \\
&\leq (1-\lambda) \|\theta^{t(n)}(t(n+m+1)) - \bar{\theta}(t(n+m+1)) +\lambda \|\theta^{t(n)}(t(n+m))-\bar{\theta}(t(n+m))\| \nonumber \\
& + \max(\lambda, 1- \lambda) \int_{t(n+m)}^{t(n+m+1)} \|\tilde{h}(\theta^{t(n)}(s),\mu(s))\|ds \nonumber \\
&\leq r(m+1,n) e^{\sum_{k=0}^{m} a(n+k) L(Y_{n+k})} + a(n+m)\left[M + C L(Y_{n+m}) e^{\sum_{k=0}^{m} a(n+k) L(Y_{n+k})}\right]\nonumber. 
\end{align}
Therefore 
\begin{align}
\rho(n,T):= \sup_{t \in [t(n),t(n) + T]} \|\theta^{t(n)}(t) - \bar{\theta}(t)\| &\leq r(\phi(n,T+1),n) e^{\sum_{k=0}^{\phi(n,T)} a(n+k) L(Y_{n+k})} \nonumber \\
                                                                      &+ a(n)\left[M + C \sup_{0\leq m \leq \phi(n,T)} L(Y_{n+m}) e^{\sum_{k=0}^{\phi(n,T)} a(n+k) L(Y_{n+k})}\right].\nonumber 
\end{align}

Now to prove the a.s. convergence of the quantity in the left hand side as $n \to \infty$, we have using Cauchy-Schwartz 
inequality:
\begin{align}
\sum_{n=1}^{\infty} E[{\rho(n,T)}^2] \leq & 2K_T \sum_{n=1}^{\infty}\left(E\left[\left(r(\phi(n,T+1),n)\right)^4\right]\right)^{1/2} + 4M^2\sum_{n=0}^{\infty}a(n)^2 + \nonumber \\
                                          & 4C^2\sum_{n=1}^{\infty} a(n)^2 E\left[\left(\sup_{0\leq m \leq \phi(n,T)} L(Y_{n+m})\right)^2 e^{2\sum_{k=0}^{\phi(n,T)} a(n+k) L(Y_{n+k})}\right],\nonumber
\end{align}
where $K_T = \sqrt{\sup_n E[e^{4\sum_{k=0}^{\phi(n,T)} a(n+k) L(Y_{n+k})}]}$ which depends only on $T$ due to \textbf{(S5)}.
Now, the third term in the R.H.S is clearly finite from the assumptions \textbf{(S4)} and \textbf{(S5)}.
Now we analyze the first term i.e. 
\begin{align}
\label{num}
\sum_{n=1}^{\infty}\left(E\left[{r(\phi(n,T+1),n)}^4\right]\right)^{1/2} \leq & 2\sqrt{2}\sum_{n=1}^{\infty}\left(E\left[\left(\sum_{k=0}^{\phi(n,T)} c_k\right)^4\right]\right)^{1/2} \nonumber \\
                                                      & + 2\sqrt{2}\sum_{n=1}^{\infty} \left(E\left[\left(\sum_{k=0}^{\phi(n,T)} a(n+k) \|M_{n+k+1}\|\right)^4\right]\right)^{1/2}.
\end{align}

Next we analyze the first term in the R.H.S of (\ref{num}) again using Cauchy-Schwartz inequality:
\begin{align}
&\sum_{n=1}^{\infty}\left(E\left[\left(\sum_{k=0}^{\phi(n,T)} c_k\right)^4\right]\right)^{1/2} \nonumber \\
&\leq 8 M^2 \sum_{n=1}^{\infty} \phi(n,T)^2 a(n)^4  \left(E\left[\left(\sup_{0 \leq k \leq \phi(n,T)} L(Y_{n+k})\right)^4\right]\right)^{1/2}+ \nonumber \\
& 8 C^2 \sum_{n=1}^{\infty} \phi(n,T)^2 a(n)^4  \left(E\left[\left(\sup_{0 \leq k \leq \phi(n,T)} L(Y_{n+k})\right)^8 e^{4\sum_{i=0}^{\phi(n,T)} a(n+i) L(Y_{n+i})}\right]\right)^{1/2}.\nonumber 
\end{align}
Therefore the the R.H.S will be finite if we can show that 
$\sum_{n=1}^{\infty} \phi(n,T)^2 a(n)^4$ is finite. 
For common step-size sequence $a(n) =\frac{1}{n}$, $\phi(n,T)= O(n)$ thus the above series converges clearly. 
One can make the series converge 
for all $a(n)= \frac{1}{n^k}$ with $\frac{1}{2} < k \leq 1$ by putting 
assumptions on higher moments in \textbf{(S4)} and \textbf{(S5)} .

In the above we have used the following inequality repeatedly for non-negative random variables $X$ and $Y$:
\begin{equation}
\sqrt{E\left[\left(X+Y\right)^{2^n}\right]} \leq 2^{\frac{2n-1}{2}}\left[\sqrt{E[X^{2^n}]} +  \sqrt{E[Y^{2^n}]}\right]\nonumber
\end{equation}
with $n \in \mathbb{N}$.

Now, 
\begin{align}
\sum_{n=1}^{\infty} \left(E\left[\left(\sum_{k=0}^{\phi(n,T)} a(n+k) \|M_{n+k+1}\|\right)^4\right]\right)^{1/2} \nonumber \\
\leq \sum_{n=1}^{\infty} a(n)^2\left(E\left[\left(\sum_{k=0} ^ {\phi(n,T)} \|M_{n+k+1}\|\right)^4\right]\right)^{1/2}\nonumber
\end{align}
which is finite under assumption \textbf{(S5)} and the fact that $a(n)$ are non-increasing.
\end{pf}
\section{Application : Off-policy temporal difference learning with linear function approximation}
\label{app}

In this section, we present an application of our results in the setting of off-policy 
temporal difference learning with linear function approximation.
In this framework, we need to estimate the value function for a target policy $\pi$ 
given the continuing evolution of the underlying MDP (with finite state and action spaces $S$ and $A$ respectively,   
specified by expected reward $r(\cdot,\cdot,\cdot)$ 
and transition probability kernel $p(\cdot|\cdot,\cdot)$) for a 
behaviour policy $\pi_b$ with $\pi \neq \pi_b$. The authors of \cite{sutton1,sutton,maeith} 
have proposed two approaches to solve the problem:
\begin{enumerate}[label=(\roman*)]
 \item  Sub-sampling: In this approach, the transitions which are relevant to deterministic target policy are kept and the rest 
of the data is discarded from the given ``on-policy'' trajectory. We use the triplet $(S,R,S')$ to represent 
(current state, reward, next state).
 Therefore one has ``off-policy'' data $(X'_n,R_n,W_n), n\geq 0$
where $E[R_{n}|X'_{n}=s,W_{n}=s'] = r(s,a,s')$, $P(W_{n} =s'|X'_{n}=s) = p(s'|s,a)$ with $\pi(s) =a$, $\pi$ 
being the target policy
 and $X'_n, n\geq 0$ is a 
random process generated by sampling the ``on-policy'' trajectory at increasing stopping times.     
 \item  Importance-weighting: In this approach, unlike sub-sampling, all the data from the given ``on-policy'' trajectory
is used. One advantage of this method is that we can allow the policy to be randomized 
in case of both behaviour and target policies unlike
the sub-sampling scenario where one can use only deterministic policy as a target policy.  
\end{enumerate}
Then they introduce gradient temporal difference 
learning algorithms (GTD) \cite{sutton1,sutton,maeith} for both the approaches.
\\ \indent
Currently, all GTD algorithms make the assumption that data is 
available in the ``off-policy'' setting i.e. of the form $(X'_n,R_n,W_n),n\geq0$ 
where $\{X'_n\}$ are i.i.d, $E[R_{n}|X'_{n}=s,W_{n}=s'] = r(s,a,s')$ and 
$P(W_{n} =s'|X'_{n}=s) = p(s'|s,a)$ with $\pi(s) =a$, $\pi$ 
being the deterministic target policy. Additionally, the distribution of $\{X'_n\}$ is 
assumed to be sampled according to the stationary distribution of the Markov chain corresponding to the behaviour policy. 
However, such data cannot be generated 
from sub-sampling given only the ``on-policy'' trajectory.   
The reason is that a Markov chain sampled at increasing stopping times cannot be i.i.d. 
In the following, we show how gradient temporal-difference learning along 
with importance weighting can be used to solve the off-policy convergence problem stated above for TD when only the 
``on-policy'' trajectory is available.
\subsection{\textit{Problem Definition}}

Suppose we are given an on-policy trajectory $(X_{n},A_n, R_{n},X_{n+1}), n\geq 0$ where  
$\{X_n\}$ is a time-homogeneous irreducible Markov chain with unique stationary distribution $\nu$ 
and generated from a behavior policy 
$\pi_b \neq \pi$. Here the quadruplet $(S,A,R,S')$ represents (current state, action, reward, next state). 
Also, assume that $\pi_b(a|s) > 0$ for all $s \in S, a \in A$. We need 
to find the solution $\theta^*$ for the following:
\begin{align}
\label{fixpoint}
0&=\sum_{s,a,s'}\nu(s)\pi(a|s)p(s'|s,a)\delta(\theta;s,a,s')\phi(s)\nonumber\\
 &=E[\rho_{X,A_n}\delta_{X,R_n,X_{n+1}}(\theta)\phi(X)]\nonumber\\ 
 &= b - A\theta,
\end{align}
where 
\begin{enumerate}[label=(\roman*)]
\item $\theta \in \mathbb{R}^d$ is the parameter for value function,
\item $\phi: S\to \mathbb{R}^d$ is a vector of state features,
\item  $X \sim \nu$,
\item  $0<\gamma < 1$ is the discount factor,
\item  $E[R_n|X_n=s,X_{n+1}=s'] = \sum_{a\in A}\pi_b(a|s)r(s,a,s')$,
\item  $P(X_{n+1}=s'|X=s) = \sum_{a\in A}\pi_b(a|s) p(s'|s,a)$,
\item  $\delta(\theta; s,a,s')= r(s,a,s') + \gamma \theta^T\phi(s') - \theta^T\phi(s)$
is the temporal difference term with expected reward,
\item  $\rho_{X,A_n} = \frac{\pi(A_n | X)}{\pi_b(A_n|X)}$,
\item  $\delta_{X,R_n,X_{n+1}}=R_n + \gamma \theta^T\phi(X_{n+1}) - \theta^T\phi(X)$ 
is the online temporal difference,
\item  $A=E[\rho_{X,A_n}\phi(X)(\phi(X) -\gamma\phi(X_{n+1}))^T]$,
\item  $b=E[\rho_{X,A_n}R_n\phi(X)]$.
\end{enumerate}
Hence the desired approximate value function under the target policy $\pi$ is $V_{\pi}^*={\theta^*}^T\phi$. 
Let $V_\theta = {\theta}^T\phi$.
It is well-known (\cite{maeith}) that $\theta^*$ satisfies the projected fixed point equation namely
\begin{equation}
V_{\theta}= \Pi_{\mathcal{G},\nu}T^{\pi}V_{\theta},\nonumber 
\end{equation}
where 
\begin{equation}
\Pi_{\mathcal{G}, \nu}\hat{V} =
\arg\min_{f \in \mathcal{G}} (\|\hat{V} - f\|_{\nu}),\nonumber 
\end{equation}
with $\mathcal{G} = \{V_{\theta} | \theta \in \mathbb{R}^d\}$
and the Bellman operator
\begin{equation}
T^{\pi}V_\theta(i) = \sum_{j \in S} \sum_{a\in A}\pi(a|i)p(j|i, a)\left[\gamma V_\theta(i) + r(i, a, j)\right]. \nonumber 
\end{equation}
Therefore to find $\theta^*$, the idea is to minimize the mean square projected 
Bellman error $J(\theta)= \|V_{\theta} - \Pi_{\mathcal{G},\nu}T^{\pi}V_{\theta}\|^2_{\nu}$ using stochastic gradient descent.
It can be shown that the expression of gradient contains product of multiple expectations. Such framework can be modelled by 
two time-scale stochastic approximation where one iterate stores the quasi-stationary estimates of some of the expectations and the 
other iterate is used for sampling. 

%
\subsection{\textit{The TDC Algorithm with importance-weighting}}

We consider the TDC (Temporal Difference with Correction) algorithm with importance-weighting 
from Sections 4.2 and 5.2 of \cite{maeith}. 
The gradient in this case can be shown to satisfy 
\begin{align}
-\frac{1}{2}\nabla J(\theta)&=E[\rho_{X,R_n}\delta_{X,R_n,X_{n+1}}(\theta)\phi(X)] - \gamma E[\rho_{X,R_n}\phi(X_{n+1})\phi(X)^T]w(\theta),\nonumber\\
w(\theta) &= E[\phi(X)\phi(X)^T]^{-1}E[\rho_{X,R_n}\delta_{X,R_n,X_{n+1}}(\theta)\phi(X)].\nonumber 
\end{align}Define $\phi_n = \phi(X_n)$, $\phi'_n = \phi(X_{n+1})$, $\delta_n(\theta) = \delta_{X_n, R_n, X_{n+1}}(\theta)$ and 
$\rho_n=\rho_{X_n,A_n}$.
Therefore the associated iterations in this algorithm are: 
\begin{equation}
\label{tdc}
\begin{split}
\theta_{n+1} &= \theta_n + a(n) \rho_n\left[\delta_{n}(\theta_n)\phi_n - \gamma \phi'_{n}\phi_n^Tw_n\right],\\
w_{n+1} &= w_n + b(n) \left[(\rho_n\delta_{n}(\theta_n) - \phi_n^Tw_n)\phi_n\right],
\end{split}
\end{equation}
\\ \indent 
with $\{a(n)\}, \{b(n)\}$ satisfying \textbf{(A4)}.
\subsection{\textit{Convergence Proof}}

\begin{thm}[Convergence of TDC with importance-weighting]
\label{th2}
Consider the iterations (\ref{tdc}) of the TDC. Assume the following:
\begin{enumerate}[label=(\roman*)]
 \item $\{a(n)\}, \{b(n)\}$ satisfy \textbf{(A4)}.
 \item $\{(X_n,R_n,X_{n+1}), n\geq0\}$ is such that $\{X_n\}$ is a time-homogeneous finite state irreducible Markov chain 
   generated from the behavior policy $\pi_b$ with unique stationary distribution $\nu$. 
$E[R_{n}|X_{n}=s,X_{n+1}=s'] = \sum_{a\in A} \pi_b(a|s)r(s,a,s')$ and $P(X_{n+1} =s'|X_{n}=s) = \sum_{a\in A}\pi_b(a|s)p(s'|s,a)$ 
where $\pi_b$ is the behaviour  policy,  
$\pi \neq \pi_b$.  
Also, $E[R_n^2 | X_n, X_{n+1}] < \infty$ for all $n$ almost surely, and
\item $C=E[\phi(X)\phi(X)^T]$ and  $A=E[\rho_{X,R_n}\phi(X)(\phi(X) -\gamma\phi(X_{n+1}))^T]$ are non-singular where $X \sim \nu$. 
\item $\pi_b(a|s) > 0$ for all $s \in S, a \in A$.
\item $\sup_n(\|\theta_n\| + \|w_n\|) < \infty$ w.p. 1. 
\end{enumerate}
Then the parameter vector $\theta_n$ converges
with probability one as $n \to \infty$ to the TD(0) solution (\ref{fixpoint}). 
\end{thm}
\begin{pf}
The iterations (\ref{tdc}) can be cast into the framework of Section \ref{defn} 
with \begin{enumerate}[label=(\roman*)]
\item $Z^{(i)}_n = X_{n-1}$,
\item $h(\theta,w,z) = E[(\rho_n(\delta_{n}(\theta)\phi_n-\gamma \phi'_{n}\phi_n^Tw))|X_{n-1}=z,\theta_n=\theta,w_n=w]$,
\item $g(\theta,w,z)=E[((\rho_n\delta_{n}(\theta) - \phi_n^Tw)\phi_n)|X_{n-1}=z,\theta_n=\theta,w_n=w]$,
\item $M^{(1)}_{n+1}=\rho_n(\delta_{n}(\theta_n)\phi_n - \gamma \phi'_{n}\phi_n^Tw_n)-E[\rho_n(\delta_{n}(\theta_n)\phi_n - \gamma \phi'_{n}\phi_n^Tw_n)|X_{n-1}, \theta_n, w_n]$,
\item $M^{(2)}_{n+1}=(\rho_n\delta_{n}(\theta_n) - \phi_n^Tw_n)\phi_n - E[(\rho_n\delta_{n}(\theta_n) - {\phi_n}^T w_n)\phi_n|X_{n-1}, \theta_n, w_n]$,
\item $\mathcal{F}_n = \sigma(\theta_m, w_m, R_{m-1}, X_{m-1},A_{m-1}, m \leq n, i = 1, 2), n \geq 0$. 
\end{enumerate}
Note that in (ii) and (iii) we can define $h$ and $g$ independent of $n$ due to time-homogeneity of $\{X_n\}$.  
\\ \indent
Now, we 
verify the assumptions \textbf{(A1)-(A7)} (mentioned in Sections \ref{defn} and \ref{assump}) for our application:
\begin{enumerate}[label=(\roman*)]
\item  \textbf{(A1)}: $Z^{(i)}_n, \forall n, i=1,2$ takes values in compact metric space as $\{X_n\}$ is a finite state Markov chain.  
\item  \textbf{(A5)}: Continuity of transition kernel follows trivially from the fact that we have a finite state MDP. 
\begin{rmk}
In fact we don't have to verify this assumption for the special case when the Markov chain is uncontrolled and has
unique stationary distribution. The reason is that in such case \textbf{(A5)} will be used only in the proof of Lemma \ref{lemma1}.
However, if the Markov chain has unique stationary distribution Lemma \ref{lemma1} trivially follows.
\end{rmk}
 \item \textbf{(A2)} \begin{enumerate}[label=(\alph*)]
\item \begin{align}
        &\|h(\theta, w,z) - h(\theta',w',z)\|\nonumber\\
        &=\|E[\rho_n(\theta-\theta')^T(\gamma \phi(X_{n+1}) - \phi(X_n))\phi(X_n) - \gamma \rho_n \phi(X_{n+1})\phi(X_n)^T(w-w')|X_{n-1}=z]\|\nonumber\\
        &\leq L(2\|\theta-\theta'\|M^2 + \|w-w'\|M^2)\nonumber,
       \end{align}
where $M=\max_{s\in S}\|\phi(s)\|$ with $S$ being the state space of the MDP and $L=\max_{(s,a)\in (S\times A)}\frac{\pi(a|s)}{\pi_b(a|s)}$. 
Hence $h$ is Lipschitz continuous in the 
first two arguments uniformly w.r.t the third. In the last inequality above, we use the 
Cauchy-Schwarz inequality. 
\item As with the case of $h$, $g$ can be shown to be Lipschitz continuous in the 
first two arguments uniformly w.r.t the third.
\item Joint continuity of $h$ and $g$ follows from (iii)(a) and (b) 
respectively as well as the finiteness of $S$.
\end{enumerate}
\item \textbf{(A3)}: Clearly, $\{M_{n+1}^{(i)}\}, i=1,2$ are martingale difference sequences w.r.t. increasing $\sigma$-fields $\mathcal{F}_n$.
Note that $E[\|M_{n+1}^{(i)}\|^2 | \mathcal{F}_n] \leq K(1 + \|\theta_n\|^2 + \|w_n\|^2)$ a.s., $n\geq 0$ since
$E[R_n^2 | X_n, X_{n+1}] < \infty$ for all $n$ almost surely and $S$ is finite.
\item \textbf{(A4)}: This follows from the conditions (i) in the statement of Theorem \ref{th2}.  
\end{enumerate}

Now, one can see that 
the faster o.d.e. becomes 
\begin{equation}
\dot{w}(t)=E[\rho_{X,A_n}\delta_{X,R_n,X_{n+1}}(\theta)\phi(X)] - E[\phi(X)\phi(X)^T]w(t).\nonumber 
\end{equation} 
Clearly, $C^{-1}E[\rho_{X,A_n}\delta_{X,R_n,X_{n+1}}(\theta)\phi(X)]$
is the globally asymptotically stable equilibrium of the o.d.e. Moreover, 
$V'(\theta,w) = \frac{1}{2} \|Cw - E[\rho_{X,A_n}\delta_{X,R_n,X_{n+1}}(\theta)\phi(X)]\|^2$ is 
continuously differentiable. Additionally, $\lambda(\theta)=C^{-1}E[\rho_{X,A_n}\delta_{X,R_n,X_{n+1}}(\theta)\phi(X)]$ 
and it is   
Lipschitz continuous in $\theta$, verifying \textbf{(A6)'}. For the slower o.d.e., the global 
attractor is $A^{-1}E[\rho_{X,A_n}R_n\phi(X)]$ verifying the additional assumption in Corollary \ref{main_col}. 
The attractor set here is a singleton. 
Also, \textbf{(A7)} is $(v)$ in the statement of Theorem \ref{th2}. 
Therefore 
the assumptions $(\mathbf{A1}) - (\mathbf{A5}), (\mathbf{A6'}), (\mathbf{A7})$ are verified. The proof would then follow
from Corollary \ref{main_col}. 
\end{pf}

\begin{rmk}
The reason for using two time-scale framework for the TDC algorithm is to make sure that the 
o.d.e's have globally asymptotically stable equilibrium. 
\end{rmk}
\begin{rmk}
Because of the fact that the gradient is a product of two expectations the scheme 
is a  ``pseudo''-gradient descent which helps to find 
the global minimum here.
\end{rmk}
\begin{rmk}
Here we assume the stability of the iterates (\ref{tdc}). 
Certain sufficient conditions have been sketched for showing  
stability of single timescale stochastic recursions with controlled Markov noise 
\cite[p.~75, Theorem 9]{borkar1}. This subsequently needs to be 
extended to the case of two time-scale recursions. 
   
Another way to ensure boundedness of the iterates is to use a projection operator. 
However, projection may introduce spurious fixed points on the boundary of the projection region and  
finding globally asymptotically stable equilibrium of a projected o.d.e. is hard. 
Therefore we do not use projection in 
our algorithm. 
\end{rmk}
\begin{rmk}
Convergence analysis for TDC with importance weighting along with 
eligibility traces cf. \cite[p.~74]{maeith} where it is called GTD($\lambda$)can be done similarly using our results. 
The main advantage is that it works for $\lambda < \frac{1}{L\gamma}$ ($\lambda\in [0,1]$ being the eligibility function) 
whereas the analysis in \cite{yu} is shown  
only for $\lambda$ very close to 1.   
\end{rmk}
\begin{rmk}
One can analyze this algorithm when the state space is infinite by imposing assumptions on $\phi$ as 
well as the target and behavior policies.
\end{rmk}

\section{Conclusion}

We presented a general framework for two time-scale stochastic approximation with controlled Markov noise.  
Moreover, using a special case  of our results, i.e., when the random process is a finite state irreducible time-homogeneous Markov 
chain (hence has a unique stationary distribution) and uncontrolled (i.e, does not depend on iterates), we provided a rigorous 
proof of convergence for off-policy temporal difference learning algorithm that is also 
extendible to eligibility traces (for a sufficiently large 
range of $\lambda$) with linear function approximation under the assumption 
that the ``on-policy'' trajectory for a behaviour policy is only available. 
This has previously not been done to our knowledge. 
\\ \indent

%
%
%

\section*{Acknowledgments.}
The authors want to thank Csaba Szepesv\'{a}ri for some useful discussion on the literature of off-policy learning. 
Our work was partly supported by the Robert Bosch Centre for Cyber-Physical Systems, Indian
Institute of Science, Bangalore.


\bibliography{mybib} 


\end{document}